\documentclass[12pt]{amsart}
\usepackage{geometry}                
\usepackage{graphicx}

\usepackage{amsmath}
\usepackage{amsfonts}
\usepackage{enumerate}
\usepackage{amssymb}
\usepackage{epstopdf}

\theoremstyle{definition}

\numberwithin{equation}{section}

\newtheorem{Thm}{Theorem}[section]
\newtheorem{Lem}[Thm]{Lemma}
\newtheorem{Pro}[Thm]{Proposition}
\newtheorem{Cor}[Thm]{Corollary}   
\newtheorem{Def}[Thm]{Definition}

\newcommand{\Diam}{\operatorname{Diam}}
\newcommand{\fll}{\operatorname{fill}}
\newcommand{\Fl}{\operatorname{Fl}}
\newcommand{\filllength}{\operatorname{filllength}}
\newcommand{\length}{\operatorname{length}}
\newcommand{\const}{\operatorname{const}}
\newcommand{\inj}{\operatorname{inj}}
\newcommand{\Vol}{\operatorname{Vol}}
\newcommand{\Area}{\operatorname{Area}}

\renewcommand{\epsilon}{\varepsilon}

\renewcommand{\rho}{\varrho}

\begin{document}

\title{Complexity of unknotting of trivial $2$-knots.}

\author{ Boris Lishak and Alexander Nabutovsky }
\maketitle

\begin{abstract}
We construct a family of trivial $2$-knots $k^i$ in $\mathbb{R}^4$ such that the maximal complexity of $2$-knots in any isotopy connecting $k^i$ with the standard unknot grows faster than a tower of exponentials of any fixed height of the complexity of $k^i$.
\par
Here we can either construct $k^i$ as smooth embeddings and measure their complexity as the ropelength (a.k.a the crumpledness) or construct PL-knots $k^i$, consider isotopies through PL knots, and measure the complexity of a PL-knot as the minimal number of flat $2$-simplices in its triangulation.
\par
These results contrast with the situation of classical knots in $\mathbb{R}^3$, where every unknot can be untied through knots of complexity that is only polynomially higher than the complexity of the initial knot.
\end{abstract}

\section{Main results.} \label{Main}

Let $k$ be a PL-unknot in $\mathbb{R}^3$ with $N$ crossings in one of its plane projections. The results of \cite{dynnikov} imply that $k$ can be isotoped to the standard unknot through PL-unknots with at most $2(N+1)^2$ crossings. 
On the other hand it was proven in \cite{nabwein3} that for each $n\geq 3$ and each {\it computable} function $f$ there exists a trivial knot $k:S^n\longrightarrow \mathbb{R}^{n+2}$ triangulated into $N$ (flat) $n$-simplices such that any isotopy between $k$ and the standard unknot that passes through $PL$-knots must pass through a knot that cannot be triangulated into less than $f(N)$ simplices.
\par
Alternatively, one can consider smooth embeddings of $S^n$ into $\mathbb{R}^{n+2}$ (or $S^{n+2}$) and measure the complexity of knots as their ropelength
(also known as crumpledness - see \cite{nabutovsky2}) that was defined as $\Vol(k)^{1\over n}\over r(k)$, where $\Vol(k)$ is the volume of $k$, and $r(k)$ denotes the injectivity radius of the normal exponential map for $k$. In other words, $r(k)$ is the supremum of all $x$ such that any two normals to $k$ of length $\leq x$ do not intersect. Informally speaking, one can think of $r(k)$ as the maximal radius of a nonself-intersecting tube centered at $k$. The same quantity $r(k)$ is also called the {\it reach} of $k$ and the thickness of $k$.
(In particular, for a $2$-knot $k$ the crumpledness is $\sqrt{\Area(k)}\over r(k)$, where $\Area(k)$ is the area of the $2$-sphere $k$.) For this measure of complexity it will still be true that if $n=1$, then there exists a polynomial upper bound for the complexity of knots in an optimal isotopy connecting an unknot with the standard unknot, and for $n>2$ the worst case complexity of knots in the optimal isotopies grows faster than any computable function.
\par
It is natural to conjecture that the results of \cite{nabwein3} for $n>2$ will also hold in the case $n=2$. Here we will prove that the complexity of untying of a trivial $2$-knot can grow faster than a tower of exponentials of any fixed height as a function of the complexity of the unknot. 

\par
\begin{Thm} \label{smoothmain} {\bf Smooth case.}
For each positive integer $k$ and arbitrarily large $N$ there exists a trivial smooth $2$-knot with crumpledness $x\geq N$ in the $4$-dimensional Euclidean space or the standard round $4$-sphere such that any smooth isotopy between this knot and the standard $2$-sphere passes through $2$-knots of crumpledness $\geq 2^{2^{\ldots 2 ^ x}}$ ($k$ times).
\end{Thm}

For a PL version of this result we will use the following notion of complexity.

\begin{Def}
Let $F$ be a PL-map between two polyhedra or two polyhedral pairs. We define the complexity $C(f)$ of $f$ as the minimal number of top-dimensional simplices in a subdivision of the domain of $f$ such that the restriction of $f$ on each simplex of the subdivision is linear.
\end{Def}

Similarly to the smooth case we will consider PL-embeddings $S^2 \to S^4$ and their isotopies $K: [0,1] \times S^2 \to S^4$. We would like to postulate the existence of unknots for which for any unknoting $K$, $\sup_{t} \{ C(K_t) \}$ grows faster than towers of exponentials. But unlike in the smooth category PL-isotopies may not be necessarily extendable to the ambient space. In fact, every PL-knot can be ``unknotted'' by an isotopy, and with only very mild increase in complexity: take the complement of a small neighbourhood of a point in $S^2$ and collapse it inside $S^4$ to a point. If we started with a non-trivial knot this isotopy would not be extendable to the ambient space while for ``complicated'' unknots an extension is possible although ``complicated''. Therefore, it makes sense to use the following definition.

\begin{Def}
Let $K: [0,1] \times S^2 \to T$ be a PL-isotopy for some polyhedron $T$ (e.g. $T = S^4$ or $\mathbb{R}^4$). We will say $K$ is locally extendable if for every $t$ there is a closed interval $I$ whose interior contains $t$ such that $K: I \times S^2 \to T$ is extendable to an isotopy $\tilde K: I \times T \to T$. We say it is $m$-effectively locally extendable if there are extensions with $C(\tilde K) \leq 2^{2^{\ldots 2 ^ {C(K_t)}}}$ ($m$ times).
\end{Def}

Local extendability implies global extendability, but the complexity of the global extension is not controlled: complexity will accumulate over time. Therefore, the following theorem is stronger than just postulating existence of unknots with fast growing complexity of global \emph{ambient} isotopy.

\begin{Thm} \label{plmain} {\bf PL case}
For arbitrary positive integer numbers $m, M, N$ there exists a trivial simplicial $2$-knot $k$ in $S^4$ or $\mathbb{R}^4$ such that $C(k)\geq N$ and for each $m$-effectively locally extendable isotopy $K: [0,1] \times S^2\longrightarrow S^4$, $K_0=k$, between this knot and the standard unknot, the complexity $C(K_t)$ at some moment $t$ is greater than $\geq 2^{2^{\ldots 2 ^ {C(k)}}}$ ($M$ times).
\end{Thm}

We will also prove a theorem of more combinatorial flavour, which does not invoke the notion of local extendability:

\begin{Thm} \label{plcomb} {\bf PL case (second version).}
For arbitrary positive integer numbers $m, M, N$ there exists a trivial $2$-knot $k: S^2 \to S^4$ (or $k: S^2 \to \mathbb{R}^4$) with complexity $x = C(k) \geq N$ such that for any sequence of knots $K_t$ ($t = 1,\ldots ,n$), where $K_1 = k$ and $K_n$ is the standard embedding, there exists $i$ such that $C(K_i) \geq 2^{2^{\ldots 2 ^ x}}$ ($M$ times), as long as for each $j\in\{1,\ldots, n-1\}$ there exists a PL-homeomorphism $\phi_j: (S^4, K_j)\longrightarrow (S^4, K_{j+1})$ with the following property: $C(\phi_j)$ does not exceed $2^{2^{\ldots 2 ^ {C(K_j)}}}$ ($m$ times).
\end{Thm}

One can use this theorem to bound from below the number of combinatorial moves needed to simplify PL $2$-knots in $\mathbb{R}^4$. For example, one can consider ``embedded'' Pachner moves: linearly embed a $3$-simplex next to the $2$-knot, so that their intersection is one, two or three, $2$-faces of the tetrahedron, 
and replace the union of these triangles with the complementary part of the boundary of the tetrahedron. The  $1$-dimensional analogue of these moves was used by Reidemeister as the definition of knot equivalence to prove his famous theorem. Theorem \ref{plcomb} immediately implies that there exist trivial $2$-knots $k$ with arbitrarily large complexity that can be unknotted only by an extremely large number of such ``embedded" Pachner moves: This number grows faster than any tower of exponentials of fixed length of $C(k)$.

\par
It is plausible that we can similarly bound below the number of PL Roseman moves required to untie our knots. Roseman moves are analogous to Reidemeister moves but in dimension one higher \cite{roseman}. And, indeed, it is possible to apply Theorem \ref{plcomb} to this case, if we can \emph{effectively} manipulate the knot diagram ($2$-dim surface with singularities in $\mathbb{R}^3$) before each Roseman move so that the Roseman move is done in some standard way of uniformly bounded complexity. It seems likely that this is possible due to the dimension and the co-dimension of the diagram being low, but we are not going to explore this in the present paper.

\par
In order to prove these theorems we first construct an infinite sequence of finite presentations of $\mathbb{Z}$. These finite presentations have certain algebraic properties that help realize them as ``visible" finite presentations of $2$-knots of complexity comparable with the total length of the corresponding finite presentations. Moreover, these finite presentations have the following additional property: In each of them there exists a trivial element of length $m$, such that one needs to apply the relations at least $2^{2^{2^{\ldots 2}}}$ times, where the height of the tower of exponentials is bounded below by $\const \ln m$ for some positive constant $\const$, in order to demonstrate that this element is, indeed, trivial. Moreover, the total length of the presentation is between $2m$ and $4m$. In the notations of the next section these presentations of $\mathbb{Z}$ are $P_{n-1}$, the trivial element is $v_n$, and $m=2^{n+1}$.
\par
Then in the smooth case we prove that the effective compactness of the set of trivial $2$-knots of bounded complexity (modulo the group of transformations of $\mathbb{R}^4$ generated by dilations and translations) implies that if these trivial $2$-knots could be ``untied" without a very large increase of complexity, then we would be able to contract any null-homotopic closed curve in the complement to the original $2$-knot to a point through closed curves that are not much longer than the original one. An analogous  assertion can also be proven in the PL-case. But the algebraic property of the finite presentations of $2$-knot groups explained in the previous paragraph implies that this is not the case.

\section{Finite presentations of $\mathbb{Z}$.} \label{Group}

The group presented by $\langle x_1, x_2\vert x_1^{x_2}=x_1^2 \rangle$ is called
the Baumslag-Solitar group. (Here and below we use the standard notation
$x^y$ for $yxy^{-1}$.) Note that for each $m$, ${x_1}^{{x_2}^m}={x_1}^{2^m}$, and 
therefore the commutator $w_m=[{x_1}^{{x_2}^m}, x_1]=e$. However, one needs to apply the relation $C \cdot 2^m$ (where $C$ is independent of $m$) times in order to demonstrate this fact in the most obvious way.
In fact, it is well-known that there is no essentially shorter
way to write $w_m$ as a product of the conjugates of the relator and its inverse. In other words, the Dehn function
of the Baumslag-Solitar group is (at least) exponential. A proof
of this fact can be found in \cite{gersten1} (and a sketch of another simpler proof
using van Kampen diagrams can be found in \cite{sapir}). 
The idea of the proof in the paper of Gersten is that the realization
complex of the finite presentation will be aspherical. (Recall that
the realization
complex of a finite presentation has one $1$-dimensional cell
corresponding to each generator of the group and one $2$-cell for each relator of the group.) So, its universal covering will be contractible, and, in particular,
will have trivial second homology group. Therefore, there will be
a unique way to fill each null homologous $1$-chain in the universal covering
by a $2$-chain. Each representation of $w_m$ as a product of conjugates of
the relator and its inverse corresponds to a filling of the lift
of the loop spelling out $w_m$ in the realization complex to
its universal covering. Therefore, it must have the same number of $2$-cells
counted with multiplicities as the number of conjugates in any presentation of $w_m$ as a product of conjugates of the relator and its inverse.
It remains only to check that $2$-cells in the universal
covering that correspond to the factors in the obvious representation of $w_m$ do not cancel, which can be done using the theory of normal forms in HNN-extensions.
The same proof implies that for each $N$ any representation
of $w_m^N$ as a product of conjugates of the relator and
its inverse has at least  $C \cdot N 2^m$ terms (because the number of cells in the $2$-chain also multiplies by $N$).
\par

One can iterate the idea used in the construction of the Baumslag-Solitar 
groups and consider the following sequence
of finite presentation of groups (see \cite{bridson1}). For each $n=1, 2,\ldots$
$$G_n=\langle x_1,\ldots, x_n\vert x_1^{x_2}=x_1^2, x_2^{x_3}=x_2^2,\ldots ,
x_{n-1}^{x_n}=x_{n-1}^2 \rangle.$$
This finite presentation has $n$ generators and $n-1$ relators of total
length $5(n-1)$. These groups are sometimes called Gromov groups.
One can show that the Dehn function of $G_n$ grows as
$2^{2^{\ldots^{2^k}}}$, where the height of the tower of exponents is $n$,
using the approach of \cite{gersten1}. (Let us denote $2^{2^{\ldots^{2^k}}}$ by $\exp_n(k)$.) The universal cover is contractible because $G_n$ is obtained as a sequence of HNN-extensions where for each extension the associated subgroups have no relations and we start with an aspherical presentation $G_1=\langle x_1 \vert \ \rangle$ (``associated subgroups'' of an HNN extension are the isomorphic subgroups mapped to each other by conjugation by the stable letter). To prove the bound on the Dehn function one can consider the following words $g_{n,k} \in G_n$ defined as $x_1^{x_2^{\ldots^{x_n^k}}}$.
It is easy to see that $g_{n,k}=x_1^{\exp_{n-1} (k)}$.
Therefore, $v_{n,k}=g_{n,k}x_1g_{n,k}^{-1}x_1^{-1}$ will be trivial. One needs to apply relations
more than $\exp_{n-1} (k)$ times to demonstrate that $v_{n,k}$ is trivial, when one proceeds in the obvious way, i.e. first eliminating $x_n$ from all copies of $g_{n,k}$, then $x_{n-1}$ etc, until we have $[x_1^{\exp_{n-1} (k)},x_1]$. As above, one can use
the asphericity of the representation $2$-complex to conclude that
the filling in its universal covering is unique on the $2$-chain level. Alternatively, one can use \cite{gersten3} or $x_i$-bands instead of asphericity. See Examples 7.2.10 from \cite{bridson1}. Note that as before with the Baumslag-Solitar group, for the area of $v_{n,k}^N$ we have a lower bound $N \exp_{n-1} (k)$, because the bound for the area of $v_{n,k}$ is homological.
Below we will choose $k=1$ to obtain $\exp_{n-1} (1)$ lower bound for the area of $v_{n,1}$.
\par
Now consider the following finite presentations $P_{n-1}=\langle G_n,t\vert tv_nt^{-1}=v_nx_n \rangle=\langle x_1,\ldots x_n, t\vert x_1^{x_2}=x_1^2,\ldots ,
x_{n-1}^{x_n}=x_{n-1}^2, tv_nt^{-1}=v_nx_n \rangle$, where $v_n=v_{n,1}=[x_1^{x_2^{\ldots x_n}}, x_1]$ are
words considered in the previous paragraph. It is easy to see that $P_n$ is a finite presentation
of $\mathbb{Z}=\langle t \rangle$ (as $v_n=e$ in the corresponding Gromov group) with the total length of relators equal to $5n-2+2^{n+2}$. The following theorem is the key technical fact in this paper.
It asserts that any way to demonstrate that, say, $x_{n+1}=e$ in $P_n$ would involve at least
$\exp_n(1)$ applications of the relations.
Equivalently, each representation of $x_{n+1}$ as the product of conjugates of the
relators and their inverses must involve at least this number of terms.
We are stating and proving this fact
using the language of van Kampen diagrams (cf. \cite{lyndonschupp}).

\begin{figure}[h]
\centering
\includegraphics[scale=0.4]{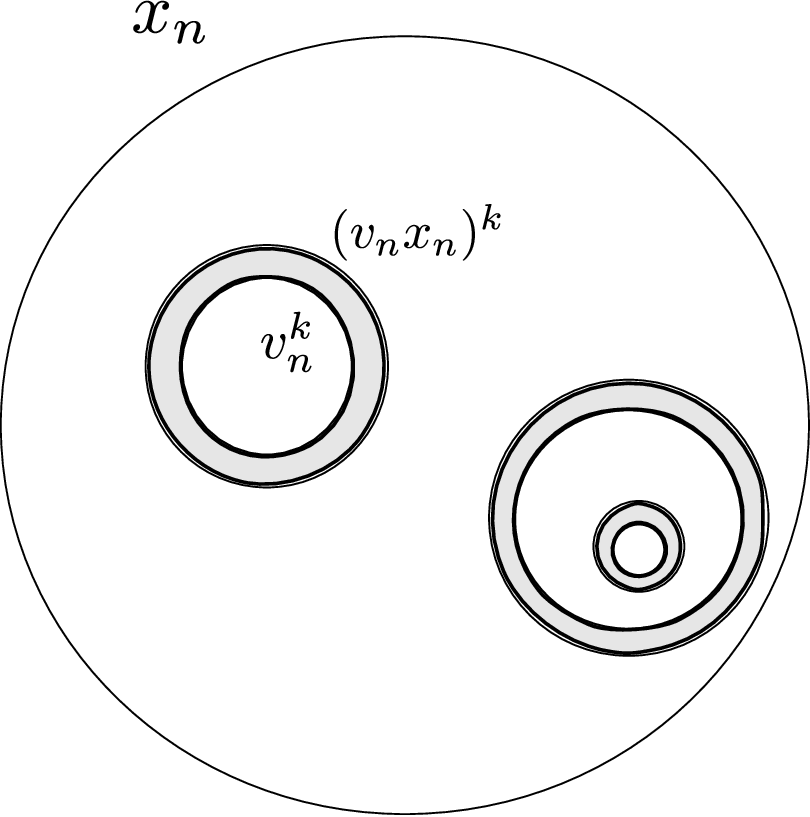}
\caption{A Van Kampen diagram for $x_n$, with some $t$-annuli marked by grey.}
\label{wn}
\end{figure}

\begin{Thm} \label{groups}
Each van Kampen diagram with the boundary $x_{n+1}$ in $P_n$ contains at least $\exp_n(1)$ cells. 
Moreover, it contains at least $\exp_n(1)$ cells corresponding
to all relators but the last one (that is, the relators not
containing the letter $t$).
\end{Thm}
\begin{proof} 
We prove this result for $P_{n-1}$ to simplify the indicies.
Consider a minimal van Kampen
diagram with $x_n$ on the outer boundary. It must contain cells corresponding
to the last relation. As there are no copies of $t$ on the boundary, these cells must form annuli ($t$-annuli) (see Figure \ref{wn}).
Consider one of the innermost $t$-annuli (that does not have any $t$-cells inside). Its inner boundary must be a non-zero power
of either $v_n$ or $v_nx_n$. The second option is impossible, as this would imply that the power of $v_nx_n$ is trivial
in $G_n$, which is false. So, we have a non-zero power of $v_n$ on the innermost boundary. The part of the van Kampen diagram
inside this innermost boundary is a van Kampen diagram in $G_n$. But we already established that any such diagram must have
size of at least $\exp_{n-1}(1)$.
\end{proof}

In the sequel we are going to use a similar result formulated for lengths rather than areas. For a finite presentation $P$ consider a sequence of words $w=F_1, ... F_n=1$, where consecutive words differ either by a free cancellation/insertion or replacing part of a relator (or its inverse) by the remaining part inverted. Therefore $w$ represents the trivial element in the group and its area is bounded from above by $n$. It is well known that we can relate the minimal length of the words in the sequence to the area of the word: if the lengths are majorized by $L$, then the number of {\it distinct} words is bounded by $const\ 2^L$, which bounds the area, or $L \geq const \log \Area$. But due to a technical complication we will need a slightly stronger (and specific for our presentations) result: 

\begin{Lem} \label{notee}
Let $x_{n+1} = F_1, ..., F_m = 1$ be a sequence of words in $P_n$ satisfying the condition above, then for some $i$ the length of the word $F_i$ not counting the letters $t, t^{-1}$ is at least $\exp_{n-2}(1)$.
\end{Lem}
\begin{proof}
Let $V$ be a Van Kampen diagram, then for any point in the diagram $p$ let $d(p)$ be the length of a shortest path in the $1$-skeleton of $V$ to the boundary of the diagram, where the length is measured as the length of the word. We denote by $\rho (V)$ the maximum of $d(p)$ over all vertices of $V$. It measures the minimal width of a neigbourhood of the boundary one needs to take to cover all of $V$.

We claim that for every van Kampen diagram $V$ with the boundary $v_n^k$ over $G_n$, $\rho(V) \geq \exp_{n-3}(1)$. We first observe that this is true for the standard diagram $D$ with the boundary $v_n$. Indeed, it is a part of the proof of the lower bound for the Dehn function of $G_n$ that this diagram injects into Cayley $2$-complex of $G_n$. (Recall, that, by definition, a Cayley $2$-complex is the universal covering of the corresponding presentation complex.)
Since $D$ sits in the Cayley complex, we know that the degree of each vertex in its image is bounded by $\const$, which allows us to bound from above the area of $r$-neighbourhoods of $\partial D$: $\length(v_n) \const^r \geq \Area_r$. Therefore, we have $\length(v_n) \const^{\rho} \geq \exp_{n-1}(1)$, and $\rho(D) \geq \exp_{n-3}(1)$.

For any diagram $V$ with the boundary $v_n^k$, its image in the Cayley complex has the boundary $v_n$ and by a homological argument has the same $2$-cells as $D$. Since $\rho$ of thus folded diagram can not decrease, we have the same lower estimate for it as for the standard disk.

Next we construct a van Kampen diagram from the sequence $F_1,...,F_m$. We can place each $F_i$ on the diagram. In it we find a subdiagram $D$ over $G_{n+1}$ with the boundary $v_{n+1}^k$ as in the proof of the theorem above. For some $i$, $F_i$ would have a subword of length $\rho(D) \geq \exp_{n-2}(1)$ inside $D$, which has no letters $t, t^{-1}$.
\end{proof}





\section{Constructing the $2$-knots.}
\subsection{Introduction}
The goal of this section is to construct the desired family of knots indexed
by $n$, and to demonstrate that their complexity (i.e. ropelength) does not grow
too fast with $n$. Our construction works as follows. We first kill the group
defined by $P_n$ by adding an extra relator $t=1$, and denote the resulting finite presentation $Q_n$. We observe that $Q_n$
is a finite presentation of the trivial group, and it can be connected with the trivial 
presentation by an explicit (and short) sequence of elementary transformations on
group presentations called Andrews-Curtis operations. (See the definition
of Andrews-Curtis operations below). The smoothed out
boundary of the $5$-dimensional $2$-handle body of each intermediate presentation will be a Riemannian $4$-sphere (see the definition of the metric below) with nice
two-sided bounds for the sectional curvature, diameter and the inverse of the
injectivity radius. The last of these spheres
will be the standard round $4$-sphere. By being somewhat more careful during
the smoothing out we can ensure that, say, $C^3$-norms of the curvature tensor
of each of these spheres will also be bounded by not too rapidly growing
function on $n$.

Andrews-Curtis operations on the finite presentations correspond to certain explicit diffeomorphisms of underlying Riemannian $4$-spheres, called handle slides and handle cancellations. These diffeomorphisms
come with isotopies connecting them with the identity diffeomorphism, but we will think of the isotopy as a family of diffeomorphisms to the target sphere with variable Riemannian metric. We will construct the unknot of a small complexity in the sphere corresponding to $Q_n$ and then pull it back to the round sphere using those diffeomorphisms. For this to be useful we need to make sure that the image of the unknot in the round sphere is of a small complexity too. We can do that by proving that 1) the diffeomorphisms are ``simple", and 2) ``simple" diffeomorphisms don't change knot complexity by much. ``Simple" for us will mean having small Lipschitz constants for the diffeomorphism and its inverse as well as for the derivatives of the diffeomorphism and its inverse.

Though 1) seems obvious if one looks at the definition of those handle slides, it's very tedious to provide all details of a rigorous proof, therefore we decided to do it differently: Each handle slide can be thought of as a smooth family of diffeomorphisms between spheres with controlled geometry. Discretizing the isotopy we can decompose
them as the composition of a large number of diffeomorphisms
so that the Lipshitz constants for both the diffeomorphism and its inverse as
well as for their derivatives are
close to $1$. Then we will need to verify the (obvious) assertion
that the number of diffeomorphisms in this discretization is ``not too
large". As this is tedious, we can use a formal argument, which does not use the details of the construction, that can be briefly described as follows (see more
details below):
Using an effective version of Gromov--Shikata theorem we can obtain a sequence of diffeomorphisms of bi-Lipschitz constant close to $1$. By a precompactness argument we can make this sequence short, and using Karcher's effective smoothing we can get a control over derivatives (see details in the next subsection).

We prove 2) now.

\begin{Lem}
Let $M^m$ be a complete Riemannian manifold, and $N^k$ its closed connected submanifold.
Assume that some real $x\geq 2$ majorizes each of the following quantities:
(a) the $C^1$-norm of the curvature tensor of $M^m$; (b) the inverse value
of the injectivity radius of $M^m$; (c) the maximal absolute value of a
principal curvature of $N^k$; (d) the distortion of $N^k$, i.e. the supremum
of the ratio of the intrinsic distance to the extrinsic distance between $x,y\in N^k$
over all pairs of distinct points of $N^k$. Then the inverse of the
injectivity radius of the normal exponential map of $N^k$, denoted $r(N^k)^{-1}$, does not exceed
$x^{c(m)}$, where $c(m)$ is a constant depending only on the dimension $m$.
\end{Lem}

\begin{proof}

First, consider the case when the ambient manifold $M^m$ is the Euclidean
space. In this case this assertion is well-known and can be found in a much 
greater generality in \cite{f}, Theorem 4.19. Here is a short proof
in the Euclidean space: we may assume that $r(N^k)^{-1}$ is greater than its maximal absolute principal curvature, because we set $c(m) \geq 1$. In that case $r(N^k)$ is small due to distant points of $N^k$ coming closer together, i.e. $2r(N^k)$ is the length of a straight line
segment $S$ perpendicular at both endpoints to the submanifold $N^k$. 
Connect these points by a minimizing geodesic $\gamma$ in the submanifold. As we have the upper bound for the distortion, the length of $\gamma$ does not exceed $2rx$. The rate of change of the unit tangent vector to the geodesic $\gamma$ is majorized by the maximal absolute principal curvature of $N^k$.
Now note that the projection of the integral of the unit tangent vector along the geodesic to 
the hyperplane perpendicular to $S$ is the zero vector. Yet at the beginning of $S$ the tangent vector is perpendicular to $S$. To compensate it must move
to the opposite unit hemisphere somewhere along the geodesic. Therefore, the geodesic cannot be too short, giving lower bound $\const(m)\over x^2$ for $r$,
and, therefore, an upper bound $x^{\const(m)}$ for $r(N^k)^{-1}$.

Now consider the general case. Without any loss of generality we can assume
that $r(N^k)$ is strictly less than the injectivity radius of $M^n$ divided by $2x$,
and that it is also less than the inverse of
the maximal absolute principal curvature of $N^k$.
In this case $2r(N^k)$ is again the length of a geodesic segment $S$ in $M^m$
connecting two points $p,q$ in $N^k$ and perpendicular to $N^k$ at each of
its endpoints. Connect $p$ and $q$ by a shortest geodesic $\gamma$
in $N^k$. The length of $\gamma$ will not exceed $2xr(N^k)$. Apply the inverse exponential map at $p$ to the closed
metric ball of radius $2xr(N^k)$ centered at $p$ in $M^m$.
Its image
will be the closed 
ball of radius $2xr(N^n)$ in the Euclidean space $T_pM^m$, the image of $p$ will be the origin, the image of $S$ will
be a straight line segment of the length $2r(N^k)$, perpendicular
to the image of $N^k$ at the origin, and the image of $\gamma$ will be
a curve $\beta$ in
the image of $N^k$. The norms of the differentials of the exponential map and the inverse
exponential map can be explicitly majorized in terms of the curvature of
 $M^m$, and the $C^2$-norms of the exponential map 
and, therefore, also of the
inverse
exponential map, can be majorized in terms of the $C^1$-norm of the curvature
tensor. Hence,
the curvature of $\beta=\exp^{-1}_p(\gamma)$ can be explicitly
estimated in terms of $x$. Now we can run the same argument, that
the integral of the unit tangent vector to $\beta$ over $\beta$ is parallel
to the straight line $\exp^{-1}(S)$, yet at the origin the tangent vector
is perpendicular
to this line, and so needs to change a lot. Juxtaposing this fact with
the observation
that $\beta$ cannot be too long we obtain the desired estimate.
\end{proof}

\begin{Cor} \label{ropelen}
Assume that $M^m$, $N^k$, $x$ satisfy the conditions of the previous lemma, $\Vol(N^k)\leq x$,
$F:M^m\longrightarrow M_1^m$ be a diffeomorphism between two Riemannian
manifolds, $x$ also majorizes the $C^1$-norm of the curvature tensor of $M_1^m$ and
the inverse of its injectivity radius. Assume that, in addition, $x$ majorizes
the Lipschitz constants of $F$, $F^{-1}$, and $DF$. Then there exists
a dimensional constant $C(m)$ such that the complexity ${\Vol^{1\over k}\over r}$
of $F(N^k)$ in $M_1^m$ does not exceed  $x^{C(m)}$. 
\end{Cor}

\begin{proof}
The diffeomorphism $F$ cannot increase $\Vol(N^k)$ by more than the factor of $x$.
The distortion cannot be changed by
more than the factor of $x^2$. The change in the maximal
absolute principal curvature can be easily and explicitly controlled in terms of the norms of the first
and the second derivatives of $F$. It remain to apply the lemma.
\end{proof}

\subsection{Construction of a family of Riemannian $4$-spheres.}
Note that if one adds one more relation to the finite presentation $P_n$, namely, $t=1$, then one obtains a finite presentation of the trivial group. Denote the resulting finite presentation of the trivial group by $Q_n$. The finite presentation $Q_n$ can be transformed to the trivial finite presentation of the trivial group by performing $O(2^n)$ elementary operations of the following types.
\begin{Def} \label{AC}
Andrews-Curtis operations:
\begin{enumerate}
    \item Inserting (or deleting) a pair of canceling generators $gg^{-1}$ or $g^{-1}g$ in any position in a generator;
    \item Replacing a relator by its inverse;
    \item Replacing a relator by its product with another relator;
    \item Cyclically permuting a relator;
    \item Adding or deleting a pair, a generator $g$ and a relator $g$, if the letter $g$ doesn't enter any other relations.
\end{enumerate}
\end{Def}

Indeed, one can use $O(2^n)$ operations involving only the relators
$tv_nt^{-1}x_n^{-1}v_n^{-1}$ and $t$ to replace these two relators by $x_n$
and $t$, then $O(1)$ operations to transform
each of the first $n-1$ relators to $x_i$ for an appropriate $i$, and finally $O(n)$ operations of type (4) removing all generators and relations to obtain the trivial presentation $\langle \ | \ \rangle$. Note the $Q_n$, the considered
trivial finite presentation, and all intermediate finite presentations
are balanced, that is the number of generators is equal to the number
of relators.
\par
For any finite presentation $P$ we can construct a smooth $4$-manifold by starting from the connected sum of several copies of $S^1\times S^3$, where $S^1$ in each copy of $S^1\times S^3$ corresponds to one of the generators, and then performing surgeries killing the relators. More precisely, we realize each relator by a  simple closed curve $\gamma$, remove the tubular neighborhood of $\gamma$, glue in a copy of $D^2\times S^2$ so that its boundary is glued to the boundary of the removed tubular neighborhood of $\gamma$ and smooth out the corners. In principle, in dimension $4$ there are two possible non-isotopic framings for a circle, but only one of those framings for each surgery will give us the alternative description of this manifold: we could start from the $2$-complex with one $0$-cell, $1$-cells corresponding to the generators and $2$-cells corresponding to relators of $P$ (i.e. the realization complex of $P$), embed it into $\mathbb{R}^5$, take the boundary of an open neighborhood of $P$ and smooth-out the corners. The resulting smooth $4$-manifold will have the fundamental group with the obvious finite presentation $P$. If $P$ is a balanced presentation of the trivial group, the second homology group is also trivial. So, it will be a homotopy $4$-sphere. Denote this manifold by $M^4(P)$. 
\par
\begin{figure}[h]
\centering
\includegraphics[scale=1.4]{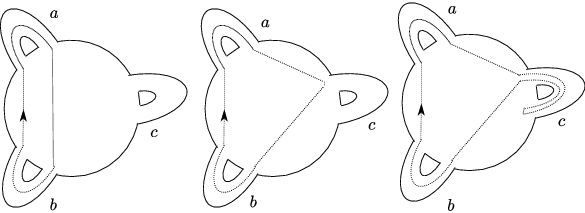}
\caption{A schematic depiction of $S^4$ with $3$ surgeries performed corresponding to generators $a,b,c$, and the isotopy of the attaching circle corresponding to move (1): the relator $ab$ is changed to $acc^{-1}b$.}
\label{finger}
\end{figure}

But it is easy to see that $M^4(Q_n)$ will be diffeomorphic to $S^4$. The reason is that if two balanced finite presentations $P$ and $Q$ are related by an elementary operation of type (3) from Definition \ref{AC}, the manifolds $M^4(P)$ and $M^4(Q)$ are diffeomorphic via a diffeomorphism that can be described as a ``handle slide", elementary operations of type (2) and (4) do not change the manifold at all, elementary operations (5) are handle cancellations, and elementary operations (1) are similar to a handle slide. More precisely, an elementary operation (1) or (3) corresponds to an isotopy of the ``attaching'' circle, a neighbourhood of which is removed to do the $2$-surgery corresponding to the relation the move is changing, see for example Figure \ref{finger} for move (1). The isotopy of the attaching circle can be extended to an isotopy of the whole $4$-manifold. After finitely many operations (we use $O(2^n)$ operations as explained above) we will end up with the standard $S^4 = M^4(\langle \ | \ \rangle)$. We will denote by $M_t$ the intermediate manifolds (the index $n$ is suppressed from this notation). Finitely many of these manifolds will correspond to the intermediate presentations of the trivial group between $Q_n$ and $\langle\ |\ \rangle$.
We will be denoting
by $\Phi_{t,s}$ the obvious diffeomorphisms between these $4$-manifolds generated by the isotopy connecting
all of them.

\begin{Def} \label{isotopy}
Let $M_t$ for $t \in [0,1]$ be the smooth family of smooth manifolds described above. In particular $M_0 = M^4(Q_n)$ and $M_1 = M^4(\langle \ | \ \rangle) = S^4$. Denote by $\Phi_{t,s}$ the diffeomorphism described above $M_t \to M_s$. Let $\phi_n = \Phi_{0,1}$.
\end{Def}

All the $M_i$ are, of course, diffeomorphic to $S^4$, but we will put different Riemannian metrics on them. We can define a Riemannian metric on $M^4(P)$ as induced from $\mathbb{R}^5$, but instead we will use a slightly different metric defined as follows. Consider a round $S^4$ of radius $1$; $n+2$ Riemannian $D^1 \times S^3$, where $D^1$ has length $1$, $S^3$ is also round and of radius $ {1 \over 1000(n+2)}$; and $n+2$ $D^2 \times S^2$, where the first factor is flat and its circumference is equal to the length of the corresponding relator and $S^2$ is round of radius ${1 \over 1000^2(n+2)L}$, where $L$ is the length of the presentation. Note that the actual Riemannian metric on these components will be slightly different, see below.
\par

We will perform surgeries on $S^4$ now as described above. First, we embed the attaching spheres of the $1$-handles ($S^0$) in $S^4$. We want to place these $2(n+2)$ points so that they are at least ${1 \over 100(n+2)}$ away from each other. This spacing allows us to embed the $3$-spheres around them. So far the resulting metric is not $C^{\infty}$: to make it smooth we need to shrink the radius of those $S^3$ in $D^1 \times S^3$ close to the endpoints of $D^1$, as if it is a continuation of the cut out $D^4$. We shrink the radius over the distance ${1 \over 10000(n+2)}$ and then smoothly restore it back to the original value of ${1 \over 1000(n+2)}$. That defines a Riemannian metric on the connected sum of $n+2$ copies of $S^1\times S^3$.

\par

Next we embed the attaching circles of the $2$-handles. We run them parallel in the boundaries of $1$-handles and then connect the end points of parallel pieces inside of $S^4$ by (almost) geodesics. We want to space out these circles so that they never approach each other closer than ${1 \over 1000(n+2)L}$. Additionally we make sure that the turns are smooth to control the curvature. We are ready now to do the surgery. We will slightly change the geometry of $D^2 \times S^2$ as follows. Firstly, we will increase the length of the circumference of $D^2$ to allow for travel inside of $S^4$, secondly, the metric is not going to be exactly the product metric: different fibres $D^2 \times {p}$ are of slightly different size, depending on the lengths of $S^1 \times {p}$ inside $S^4$, and lastly we will smooth out this metric to make it $C^{\infty}$ similarly to how we did for $0$-surgeries in the previous paragraph. Now we have a Riemannian metric on $M(P)$ a finite presentation $P$.
\par

We can similarly define Riemannian metric on such manifolds during the moves. As we move the attaching circle we continue to make sure that it bends smoothly and there is sufficient spacing around it. Since circles cannot knot in dimension $4$, the only obstruction is the available space, but since the size of the buffer zone decreases with the length of the presentation we have plenty of space. We scale the size of the $2$-handle according to the length of the attaching circle. This defines Riemannian metrics on $M_i$. 
\par

There are some ambiguities in defining this metric. First, we didn't specify where exactly we place the attaching spheres $S^0$ inside $S^4$ and which points of the second factor of $D^1 \times S^3$ inside correspond to the attaching spheres $S^1$. Secondly, we didn't specify how the attaching spheres $S^1$ with their neighbourhoods avoid each other inside $S^4$. The diffeomorphism class does not depend on these choices for the same reason as handle slides discussed above induce ambient isotopies. But also the defined metrics are close to each other in the following sense.

\par
As we smoothly vary the attaching circle and remove it together with a neighbourhood from the manifold, the resulting submanifolds are continuously close to each other in the Gromov--Hausdorff metric since they are subspaces of the same space. Because we attach a $2$-handle with a continuously scaled metric, as the circle gets shorter or longer, we can obtain a sequence of Riemannian metrics between the initial manifold and the one after the move with arbitrarily short distances, in Gromov--Hausdorff metric, between consecutive metrics. We will collect the above discussion into a lemma, but first we make the following definition:

\begin{Def} Let $f$ and $g$ be two positive valued functions
defined on a closed unbounded subset $D$ of $[0,\infty)$. We say that they have
{\it similar growth} if there exists $N$ such that $f(x)< \exp_N( g(\{ \exp_N ( x) \}_D ))$ and $g(x)< \exp_N( f(\{ \exp_N ( x) \}_D ))$. Here $\{y\}_D$ means $\min\{x \in D | x \geq y \}$. Increasing functions that do {\it not} have similar growth with at least one function $f$ such that $f(x) \leq x$ (restricted to their domain) are called {\it rapidly growing} functions. An increasing function that is not rapidly growing is called {\it reasonably growing}.
\end{Def}

Now the previous discussion easily implies that:

\begin{Lem}\label{smoothpresent}
For any $\eta>0$ there is a finite sequence $t_i$ between $0$ and $1$ such that $M_{t_i}$ and $M_{t_{i+1}}$ differ from each other by no more than $\eta$ in the Gromov--Hausdorff metric and for any $M_t$ the diameter $\Diam$, the volume $\Vol$, sectional curvatures $K$, the injectivity radius $\inj$ have the following bounds: $1<\Vol<V$, $\Diam <D$, $|K| < F$, $\inj > {1 \over G}$, where $V, D, F, G$ are reasonably growing functions of $n$ independent of $\eta$. Moreover, the $C^1$-norm of the curvature tensor of $M_{t_i}$ will not exceed $F(n)$ for all $t_i$.
\end{Lem} 

\begin{proof}

We can get control over $C^1$-norm by taking additional care in how smooth are the attaching circles and how we smooth out the metric after surgeries. (This is not surprising, as any closed Riemannian manifold with $\vert K\vert\leq 1$ can be effectively smoothed out using the Ricci flow to gain a control over any number of derivatives of the curvature tensor, while staying arbitrarily $C^0$-close to the original manifold - cf. \cite{bmr} or \cite{nabwein1} for more details.)
\end{proof}

{\bf Remark.} One can see that consecutive manifolds in the sequence from this lemma are bi-Lipschitz homeomorphic by a diffeomorphism $\Phi$ such that the Lipschitz constants for $\Phi$ and its inverse do not exceed $1 +\eta$ directly from construction, but it also formally follows from one of the (constructive) proofs (cf. \cite{bmn}, which uses \cite{k}) of the Gromov--Shikata theorem \cite{gromov}, Theorem 8.19. These proofs imply that sufficiently Gromov--Hausdorff close manifolds satisfying upper bounds for $\vert K\vert$, $\Diam$ and a lower bound for $\inj$ are bi-Lipschitz homeomorphic by means of a diffeomorphism,
and imply the concrete estimates on the Lipschitz distance in terms of the Gromov--Hausdorff distance. We are going to use this fact in the next subsection.

\subsection{Construction of the $2$-knots.}
Now we are going to construct the $2$-(un)knots. Consider the last relator, $t$, and its attaching circle $S^1$, also denoted $t$, that we can consider to be in $M(P_n)$. The boundary of the attached handle is diffeomorphic to  $S^2\times D^2$. Consider the $2$-sphere $S_n=S^2\times \{c\}$, where $c$ is the center of $D^2$. We claim that $S_n$ represents the trivial $2$-knot. Indeed, we can use operations 1)-5) on $P_n$ to reduce it to $\langle t | \ \rangle$, which induces an isotopy of $M(Q_n)$ to $M(\langle t | t \rangle)$ which fixes $S_n$ (because we didn't use the relation $t$ in our reduction), but $S_n$ is clearly trivial in $M(\langle t | t \rangle)$. This shows that the standard pair $(S^4, S^2)$ is diffeomorphic to $(M(Q_n), S_n)$, which is diffeomorphic to $(S^4, \phi_n(S_n))$, i.e. $\phi_n(S_n)$ is a family of unknots.

\par

Next we want to make sure that $\phi_n(S_n)$ has reasonable complexity. 
Corollary \ref{ropelen} implies that this boils down to verifying that $\phi_n$ can be decomposed as a composition of a reasonably growing number of
diffeomorphisms with reasonably growing $C^2$-norms and $C^1$-norms of the inverses.

This could be seen from a description we will give later of all the intermediate knots that will have nicely controlled (and clearly reasonably growing) complexities, but to turn it into a complete proof would require going through a lot of tedious details. Therefore, as mentioned in the introductory subsection, we will use the previous lemma and give some a priori bounds that are reasonably growing functions of $n$ in a way that does not appeal to the geometry of the constructed $2$-knots. The above discussion implies
that a reasonably growing in terms of $n$ upper bound for $\eta$-jumps between $M(Q_n)$ and the standard sphere
$S^4$ in the Gromov--Hausdorff metric would be sufficient for our purposes.
Here we can assume that $\eta^{-1}$ is an explicit  reasonably growing function of $n$ (one can take $\eta^{-1}=100G(n)$, where $G$ is the same as in Lemma 3.6), and all jumps must be through Riemannian manifolds 
with upper bounds for $\vert K\vert$, $\Diam$ and $\inj^{-1}$.
Now recall the well-known explicit exponential in a power of $\exp(D)\over\epsilon$ 
upper bound for the $\epsilon$-entropy of spaces of isometry
classes of Riemannian metrics satisfying a constant negative lower bound for the Ricci curvature and the upper bound $D$
for the diameter. This bound easily implies the existence of a reasonably growing in $n$ upper bound
for the $\eta/2$-entropy of the space of isometry classes of Riemannian metrics on $S^4$ with the upper bounds
for $\vert K\vert$, $\Diam$ and $\inj^{-1}$ as in the previous lemma. (By definition, the $\epsilon$-entropy of a compact metric space
is the minimal number of metric balls of radius $\epsilon$ required to cover this space, or, equivalently, the number of points in a minimal $\epsilon$-net in
the space.)

Now the key observation is that
given a sequence of $\eta$-jumps as in the previous lemma, we can replace it by a {\it shortest} sequence of $\eta$-jumps
satisfying the same condition but passing through points in a minimal $\eta/2$-net.
The length of the shortest sequence will be automatically bounded in terms
of $\eta/2$-entropy of the considered space of isometry classes of Riemannian metrics on $S^4$
satisfying the bounds on $\vert K\vert, \Diam, \inj^{-1}$ - and, therefore, is reasonably growing. (Indeed, the shortest sequence of jumps can pass through each point in the net at most once.)

Now, as in the remark after Lemma \ref{smoothpresent} every such step yields a diffeomorphism with small Lipschitz constants. We will abuse notation by redefining $\phi_n$ as the composition of these diffeomorphisms.

Finally, we can use the fact that the Riemannian $4$-spheres have nicely controlled geometry, and to use (constructive) smoothing with mollifiers from \cite{k} to replace any diffeomorphism $F$ with nicely controlled $C^0$-norms of $DF$ and $DF^{-1}$ by a close diffeomorphism, where one has, in addition, a nice and explicit control over the second derivatives of $F$ and its inverse. Applying Corollary \ref{ropelen} we get the next lemma.

\begin{Lem}\label{knots}
The knots $\phi_n(S_n)$ have reasonably growing complexity. Additionally, the Lipschitz constants for $\phi_n$ and its inverse are reasonably growing functions of $n$.
\end{Lem}

\subsection{A diagram of the $2$-knots.}
Now we are going to give a picture of these knots, though we are not going to rely on it in our proofs. A \emph{superslice} knot is a slice $1$-knot with the property that taking the double of $D^4$ with a slice disk in it produces the trivial $2$-knot. We are going to represent our $2$-unknot by such a superslice knot. See Figure 2 of \cite{meier} for the $2$-unknot similarly represented.

\begin{figure}[h]
\centering
\includegraphics[scale=0.7]{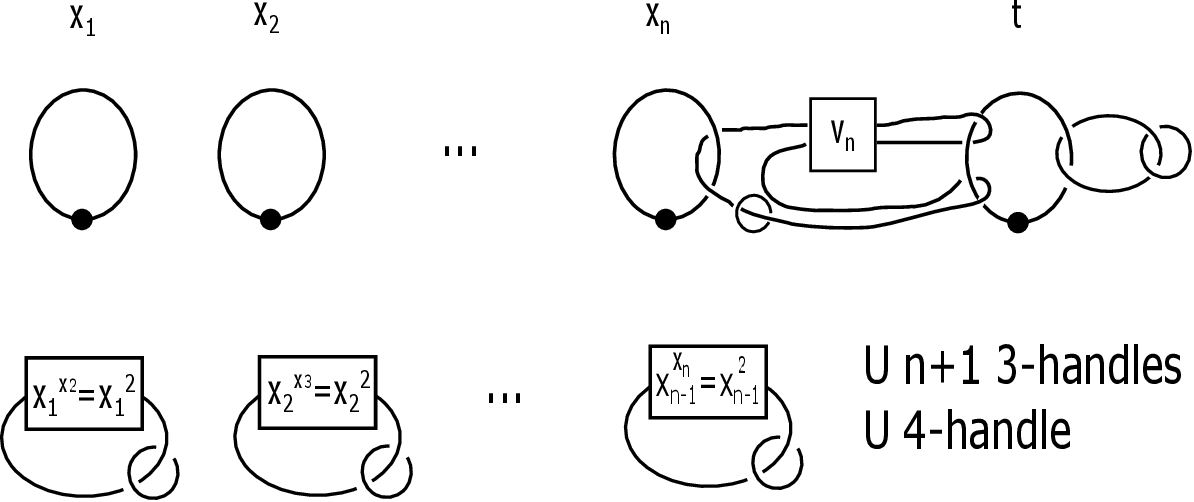}
\caption{The natural handle decomposition of $M^4(Q_n)$. All framings are $0$.}
\label{kirby}
\end{figure}

We start by giving a natural handle decomposition of $M^4(Q_n)$. This sphere is the boundary of the $5$-dimensional manifold consisting of $1$ $0$-handle, $n+1$ $1$-handles and $n+1$ $2$-handles. Equivalently, this sphere is the double of the $4$-manifold with the same handle decomposition. By looking at the second copy of the $4$-manifold upside down we see that this double has $1$ $0$-handle, $n+1$ $1$-handles, $2n+2$ $2$-handles, $n+1$ $3$-handles, and $1$ $4$-handle (see Example 4.6.3 in \cite{stips}). Every original $2$ handle gets its copy as a $0$-framed meridian. In fact, all framings are $0$ (by Proposition 5.7.1 in \cite{stips} the framings are even as the manifold is spin, and we can make them zero using the meridian). The knotting of the attaching circles is unimportant because meridian handles allow us to switch any crossing - this is expected since this manifold comes from $5$-dimensional surgery. Below we give a Kirby diagram with Akbulut notation for $1$-handles as dotted circles. By a path entering a rectangle with a word written in it we mean that this path travels through $1$-handles according to this word, see Figure \ref{kirby}.

\begin{figure}[h]
\centering
\includegraphics[scale=0.7]{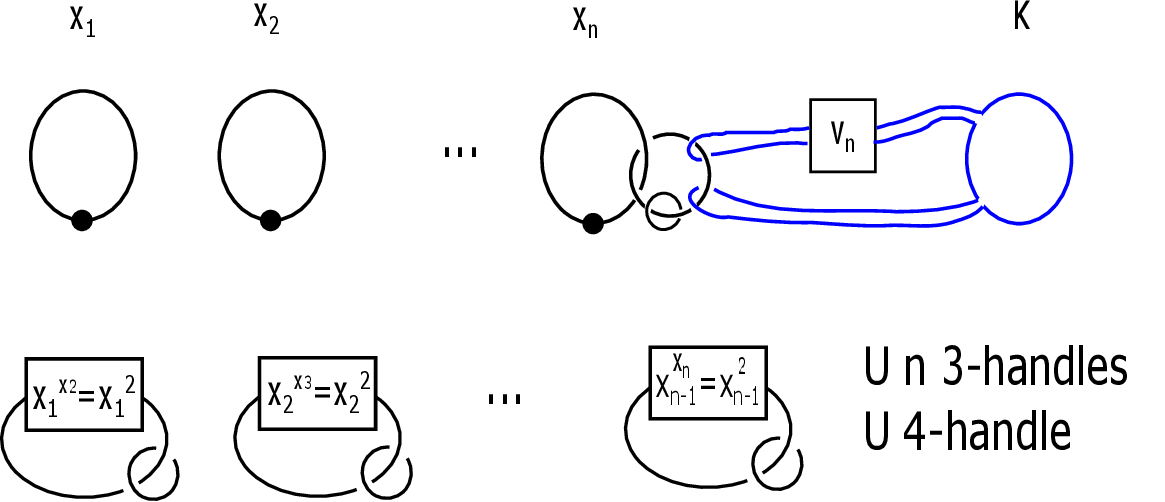}
\caption{We simplified the last relator to $x_n$, which isotoped $S_n$ to the sphere intersecting this diagram in the blue knot $K$.}
\label{knot1}
\end{figure}

We can cancel the relator $t$ with the generator $t$ (while the meridian cancels with one of the $3$-handles). Then we can see the intersection of $S_n$ with the diagram $S^3$ (call it $K$) as the same circle which used to denote the generator $t$. To see $\phi_n(S_n)$ we will continue to simplify this handle decomposition, starting with the isotopy of the last relator, see Figure \ref{knot1}.

\begin{figure}[h]
\centering
\includegraphics[scale=0.7]{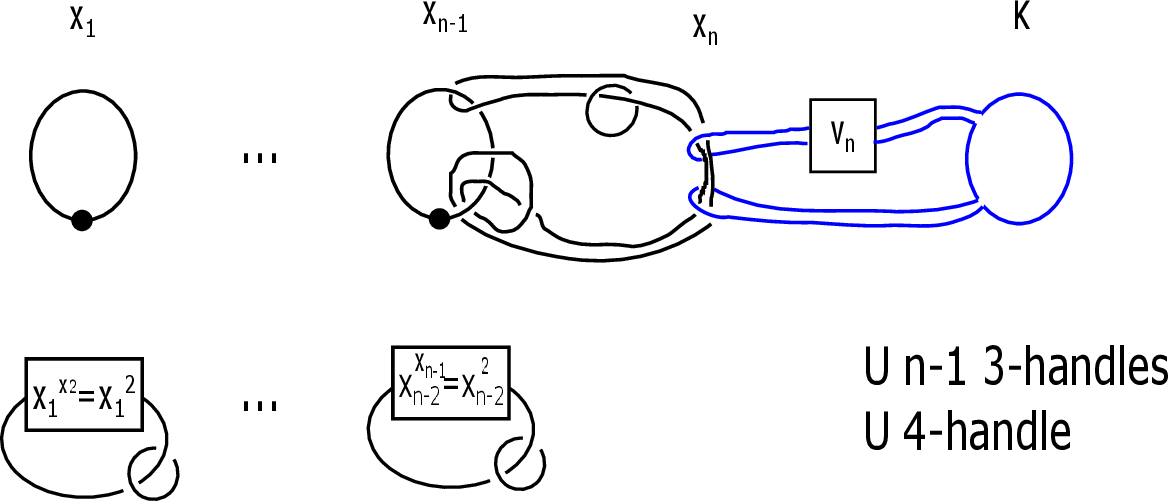}
\caption{We draw the relator $x_{n-1}^{x_n} = x_{n-1}^2$ explicitly.}
\label{knot2}
\end{figure}

Next we can handle slide every other $2$-handle off $x_n$ and cancel it with the $1$-handle. Note that after the handle slides, the two loops of $K$ now play the role of the generator $x_n$, and that's what boxes now mean: if there is the $x_n$, go through both loops with opposite orientation, see Figure \ref{knot2}.

\begin{figure}[h]
\centering
\includegraphics[scale=0.7]{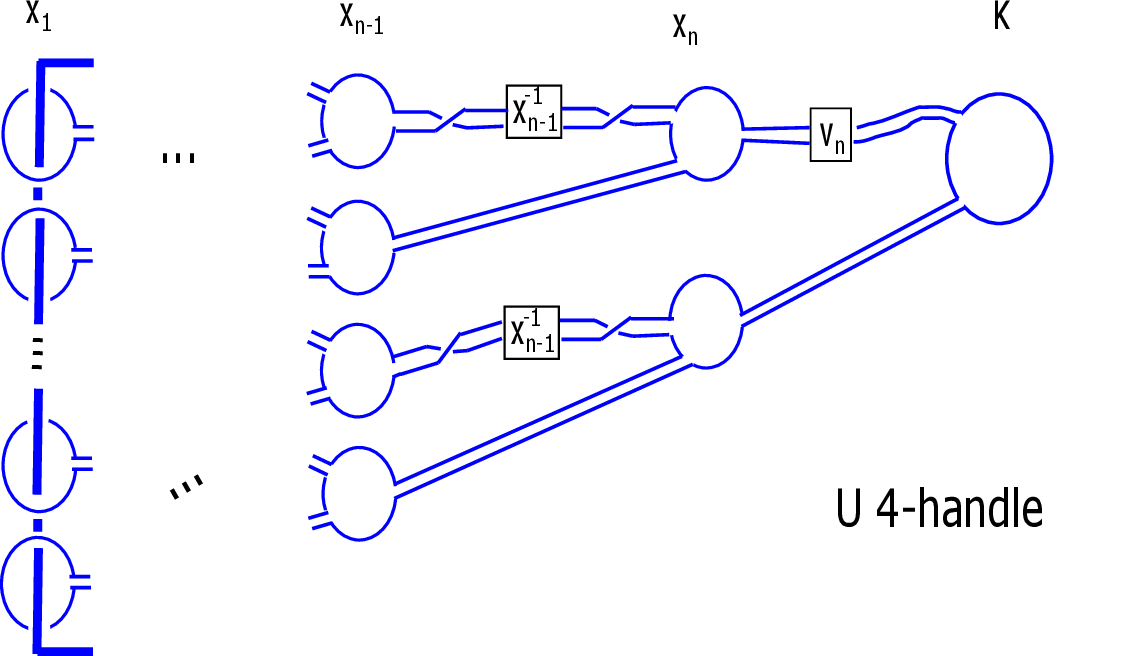}
\caption{The intersection of $\phi_n(S_n)$ with the meridian sphere $S^3$. The thick blue line denotes a bunch of $2^n$ pairs of blue strands.}
\label{knot3}
\end{figure}

We continue like this until we cancel all handles to get a picture for $\phi_n(S_n)$, see Figure \ref{knot3}. Note that double strands on this diagram represent small tubes that can pass through each other in $4$ dimensions preserving the $2$-knot, but they can't pass through single strands (See Figure 1 of \cite{meier}). For double strands that can pass through each other, this knot emulates relations in the corresponding group. This gives a very long isotopy simplifying $\phi_n(S_n)$ by simplifying $v_n$. The purpose of the next section is to prove that one can't do it much faster.

One can also perform a stereographic projection from a point on $S^4$ far from $\phi_n(S_n)$ and obtain desired $2$-knots in $\mathbb{R}^4$.

\section{Filling functions.}

Now we are going to use a concept from \cite{gromov}:

\begin{Def}

For each Riemannian manifold, or, more generally, length space $X$ of finite diameter we define its {\it filling length} function $\Fl_x(X)$ as follows. 
For a closed contractible curve $\gamma$ let $H(\gamma)$ be the set of all homotopies contracting $\gamma$ to a point. For each $h\in H(\gamma)$ let $\fll(h)$ denote the maximal length of a closed curve in $h$. Then define $\filllength_X(\gamma)$ as the $\inf_{h\in H(\gamma)} \fll (h)$. Finally, we define $\Fl_x(X)$ as the supremum of $\filllength_X(\gamma)$ over the set of all closed contractible piecewise smooth curves $\gamma$ of length $\leq x$, and {\it the filling constant of} $X$ $\Fl(X) = \sup_x \ {\Fl_x(X) \over x}$.
\end{Def}

For $M$ a compact Riemannian manifold, $\Fl_x(M)$ behaves essentially like the Dehn function of $\pi_1(M)$, which we will see in the proof of the next lemma. We are going to use $\Fl$ for knot complements with cyclic fundamental group. This group has linearly growing Dehn function which implies the finitness of $\Fl$.

\begin{Lem}
For a compact Riemannian manifold $M$ with or without boundary $\Fl_x(M)$ is finite for all $x$. If $\pi_1(M) = \mathbb{Z}$, then the filling constant of $\Fl(M)$ is also finite.
\end{Lem}

\begin{proof}
To prove the first assertion, observe that if the length of $\gamma$ is very small, then $\gamma$ can be contracted to a point without
any length increase. Now the first assertion follows from the precompactness of the set of all piecewise
closed curves of length $\leq x$ parametrized by the arclength.

To prove the second assertion
let $D$ be the diameter of $M$. Consider a fine triangulation $T$ of $M$. We can use the $1$-skeleton of $T$ and any spanning tree $t$ in this $1$- skeleton to define a set of generators $g_1,\ldots , g_n$ of $\pi_1(M)$, where $n=n(M,T,t)$. Some of the two dimensional simplices of $T$
correspond to (all) relators in a finite presentation of $\pi_1(M)$. Given a contractible loop $\gamma$,
one can first homotop it to a geodesic loop via a length non-increasing homotopy, and then to
a simplicial loop in the $1$-skeleton of $T$ via a homotopy that increases length by a factor that depends only on $M$ and $T$ but not on
$\gamma$. Now it can be represented as a word in $g_i$ and their inverses of length bounded by $\const_1(M,T,t)\length(\gamma)$.
This word represents the trivial element in $\pi_1(M)=\mathbb{Z}$, and it can be reduced to the trivial word applying less than $\const_2(M,T,t)\const_1(M,T,t)\length(\gamma)$ relations (as $\mathbb{Z}$ has a linear Dehn function). Each application of a relation increases the length of the curve by a very small summand (the length
of which is bounded by $\const_3(M,T,t)$). Taking the supremum over all contractible loops of length $\leq x$ we obtain
an upper bound on $\Fl_x(M)$ of the form $\const(M,T,t)x$. Choosing $T$ and then $t$ in any canonical way, and dividing both sides by $x$, we see that $\Fl(M)\leq \const(M)$.

\end{proof}

Now we plan to use these filling functions similarly to how it was done in \cite{nabutovsky3}. The main idea is that they behave in a similar way to their algebraic counterparts that measure how difficult it is to see that all trivial words in a ``visible" finite presentation of the fundamental group of $M$ are, indeed, trivial (or, more concretely, the maximum over all trivial words of a given length of the minimal area of a van Kampen diagram for the considered word). More specifically, we would like to introduce $N(S_n)$ defined as the tubular neighbourhood
of $S_n$ in $M(Q_n)$ of radius ${1\over 10B(n)}$, where $B(n)$ is a reasonably growing function, and ${1\over B(n)}$ is the common lower bound of $r$
for all $2$-knots in the sequence connecting $S_n$  in $M(Q_n)$
and $\phi_n(S_n)$ in the standard round $S^4$. Similarly we define $N(\phi_n(S_n))$ as the open tubular neighbourhood of $\phi_n(S_n)$
in $S^4$ of radius ${1\over 10B(n)}$.
Then we consider the
complements to $N(S_n)$ in $M(Q_n)$, and $\phi_n(N(S_n))$ in the standard $S^4$ and want to establish that:
\begin{enumerate}
    \item If knots $\phi_n$ in the standard $S^4$
can be untied through $2$-knots of 
a reasonably growing complexity, then
$\Fl(S^4\setminus N(\phi_n(S_n)))$ is reasonably growing;
    \item the values of $\Fl$ for these two $2$-knot complements (in $S^4$ and $M(Q_n)$) are similarly growing functions of $n$; and
    \item $\Fl$ for the complements of $N(S_n)$ in $M(Q_n)$ is a not reasonably growing function of $n$, as its growth is closely
    related with the growth of the area of van Kampen diagrams required to demonstrate that the groups of these $2$-knots are trivial (see Theorem \ref{groups}).

\end{enumerate}

Taken together these three facts would establish that the constructed $2$-knots can be untied only through $2$-knots of a very high complexity. The fact that we removed
not merely $\phi_n(S_n)$ but its open neighbourhood ending up with a compact
space is helpful not only because so far
$\Fl(X)$ was defined only for compact spaces
but also in a specific proof of assertion 1
given in the next proposition. On the other
hand Lemma \ref{knots} asserting 
among other things that the
diffeomorphism between $(M(Q_n), S_n)$ and $(S^4,\phi_n(S_n))$ and its inverse have
reasonably growing Lipschitz constants (regarded as functions of $n$) would immediately
imply that $\Fl(M(Q_n)\setminus S_n)$ and
$\Fl(S^4\setminus \phi_n(S_n))$ are similarly growing providing we would make sense of $\Fl$
in the non-compact case. On the contrary, it is
not that easy to compare the filling constant
of the complements of the tubular neighbourhoods of $S_n$ and $\phi_n(S_n)$
as the image of the tubular neighbourhood
under $\phi_n$ is not the same as the tubular
neighborhood of the image - and it is not
so easy to relate the two in our case. To circumvent this difficulty we prove that
$\Fl$ of the complement of a $2$-sphere $S$
in a $4$-sphere $\Sigma$ is finite and can
be majorized in terms of $\Fl$ of the 
complement to the tubular neighbourhood of
$S$ of a reasonably small thickness. As
an immediate corollary, we will have assertion 1 for
$\Fl$ of the complements to the $2$-spheres (in addition to just established assertion 2).
Finally, we will establish assertion 3
directly for $\Fl$ of the complement of $S_n$
rather than $N(S_n)$. Again, together
these three assertions for the complements of the spheres will imply our main
theorem.

To summarize: we prove Theorem \ref{smoothmain} by first obtaining a sequence of knots $\phi_n(S_n)$ from Lemma \ref{knots}. By Proposition \ref{dehn}, $\Fl(M(Q_n)\setminus S_n)$ is not reasonably growing. Then by Proposition \ref{standard} $\Fl(S^4 \setminus \phi_n(S_n))$ is not reasonably growing, and from Lemma \ref{complement} it follows that $\Fl(S^4 \setminus N( \phi_n(S_n)))$ is also not reasonably growing. Finally, by Proposition \ref{comfill}, this sequence of knots can not be easily untied.

Now we are going to present a proof of the next proposition modulo some technical lemmas that are proved right after. We are going to use $C(k)$ for the square of crumpledness of $k$, that is $C(k) = {\Area(k) \over r^2(k)}$.

\begin{Pro} \label{comfill}
Let $k^{(i)}$ be a sequence of 2-unknots in the round sphere $S^4$ of radius $1$ with ${1\over r(k^{(i)})}\leq B(i)$ for some reasonaby growing function $B$. Assume that $k^{(i)}$ can be untied through $2$-knots of complexity bounded above by a reasonably growing function of $i$. Then $\Fl(S^4 \setminus N_{1\over 10B(i)}(k^{(i)}))$ is a reasonably growing function of $i$.
\end{Pro}

\par\noindent
{\bf Remark.}
The assumption that inverses of $r(k^{(i)})$ are bounded
above by a reasonably growing function of $i$ is unnecessary here,
if instead of $N_{1\over 10B(i)}(k^{(i)})$ one would use, say, $N_{r(k^{(i)})\over 10}(k^{(i)})$.

Indeed, if $r(k^{(i)})$ is very small, it is easy to see that $k^{(i)}$
is in a metric ball of a very small radius centered at some point
$p\in S^4$, and rescaling this ball to a much larger but still small
metric ball we obtain a new knot with reasonably large $r$ and comparable complexity $C$ that can be connected with $k^{(i)}$ 
by an isotopy that passes through intermediate rescalings of
$k^{(i)}$. It is easy to see that $\Fl$ of the complements of the appropriate tubular neighbourhoods of these rescalings are comparable to $\Fl(S^4 \setminus N_{1\over 10B(i)}(k^{(i)}))$. (By a rescaling here and below we
mean a map of a smaller metric ball $B_{r_1}(p)$ to a larger
concentric ball $B_{r_2}(p)$ obtained as a composition of the inverse
exponential map to $T_p\ S^4$, rescaling by the factor of ${r_2\over r_1}$, and the exponential map $\exp_p$. The center of the rescaling should not be too close to the knot.)

Yet we do not need this stronger version of the lemma.

\begin{proof}
We suppose that $k^{(i)}$ can be untied through knots $K_t,\ t\in[0,1]$ (index $i$ suppressed) of complexity bounded by a reasonably growing function of $i$. By Lemma \ref{compactc} without loss of generality we can assume that the areas of the knots during the isotopy are between $1 \over g(i)$ and $g(i)$, where $g$ is a reasonably growing function, and the injectivity radius of the normal exponential map during a contracting isotopy is bounded below by ${1\over F(i)}$, where $F$ is also a reasonably growing function. From Generalized Gauss Lemma (see, for example, \cite{tubes}) it follows that equidistant to $K_t$ surfaces, $\partial N_{\epsilon}(K_t)$,  have radius of the normal exponential map bounded from below by ${1 \over F(i)} - \epsilon$.
\par

The idea of the proof of this proposition is to discretize the isotopy $K_t$ in a number of ``short'' steps, so that each step does not change the function $\Fl$ by much. The number of steps will be kept small using the fact that the space of knots of bounded complexity is ``small''. 

Observe that the isotopy between the standard unknot and the given unknot can be replaced by a sequence of ``jumps" of ``length" $\leq {1\over 100000F(i)}$, where the number of jumps is bounded by a reasonably growing function of $i$ for the following reason (``length'' here is the Gromov--Hausdorff distance): The number of these jumps is bounded by the number of points in an ${1\over 100000F(i)}$-net in the space of compact Riemannian manifolds $S^4 \setminus \partial N(K_t)$ plus two. Indeed, we can take every ``jump'' except the first and the last one to be between the points in the net. The number of points in such a net is a reasonably growing function of $i$ by Lemma \ref{precompact}.

\par
Now we can start at the standard unknot, ``jump" back to the given $2$-knot and observe that $\Fl(S^4 \setminus N(k^{(i)}))$ does not exceed $\const^{\#\ jumps}$ by Lemma \ref{step}, which is bounded by a reasonably growing function of $i$.
\end{proof}

Note that one can prove the result for $2$-knots in $\mathbb{R}^4$ instead of $S^4$ by scaling the isotopy $k_t$ to a very small neighbourhood of the origin and then placing this isotopy on $S^4$ without much change in the complexity. Now we prove the technical lemmas.

\begin{Lem} \label{radial}
Let $M$ be a smooth closed submanifold of $S^4$ with the standard round metric and $r \leq r(M)$. Let $f$ be the radial projection $N_{r \over 10}(M) \to M$. Then the Lipschitz constant of $f$ does not exceed $10 \over 9$. Similarly, the Lipschitz constants of the radial projection diffeomorphism $\pi:\partial N_{\delta}(M) \to \partial N_{2 \delta}(M)$ and its inverse are bounded by $2$, if $ {r \over 10} < \delta < {r \over 10} + {r \over 100}$.
\end{Lem}

\begin{proof}

The idea is that $f$ can not stretch vectors by a lot because we are far away from where the normal exponential map becomes singular. First, consider the flat $R^4$ instead of $S^4$. We would like to prove that for two infinitesimally close points $x,y\in M$
and two unit normals to $M$, $n_1$ at $x$ and $n_2$ at $y$, the distance between $x+{r\over 10}n_1$ and $y+{r\over 10}n_2$
cannot be too small in comparison with the distance between $x$ and $y$. As $x$ and $y$ are infinitesimally close,
we are studying the distance between two straight lines $x+tn_1$ and $x+tn_2$, and up to the first order can
approximate this distance by the distance between points $x+tn_1$ and $y+tn_2$, where we care about $t={r\over 10}$.
It is obvious that if we want to minimize this distance, we need to choose $n_2$ so that the projection of $n_2-n_1$ to $x-y$ is
maximal, i.e. to choose $n_2$ in the plane spanned by $n_1$ and $y-x$. In this case, the best chance to have a small distance
between $y+{r\over 10}n_2$ and $x+{r\over 10}n_1$ is when these lines are intersecting as close as possible to $x$. Yet they cannot
intersect at a distance less than $r$. Looking at the similar triangles we see that the Lipschitz constant for $f^{-1}$ is not less
than ${9\over 10}$ , and, therefore, the Lipschitz constant for $f$ does not exceed ${10\over 9}$.
\par

In the case of the ambient round $S^4$ we are studying the distance between infinitesimally close geodesics (great circles),
and similarly we can assume that they belong to the same equatorial $S^2$. As the result we will be studying the change of length of a Jacobi
fields perpendicular to one of these geodesics along the geodesic, that can vanish at the distance not less than $r$ from the starting point. Denote this distance by $r+c$, where $c\geq 0$. Then the quantity that we are interested will be ${\vert J(r+c)\vert\over \vert J({9\over 10}r+c)\vert}$.  Using the explicit
formulae for the Jacobi fields on the unit sphere, this quantity can be majorized by ${\sin(r+c) \over \sin ({9 r \over 10} +c )}$.
This quantity is minimal, when $c=0$, is equal to ${\sin r\over sin ({9\over 10}r)}$, and
does not exceed ${1\over 9}$.

\par
One can compute the upper bound for the Lipschitz constant of $\pi^{-1}$ in exactly the same fashion.
To compute the Lipschitz constant for $\pi$, we will need to majorize $\vert J(2\delta)\vert\over \vert J(\delta)\vert$, where the
Jacobi field $J$ vanishes at $0$. Therefore we will need to majorize ${\sin (2\delta)\over\sin (\delta)}=2\cos \delta<2.$
\end{proof}

\begin{Lem} \label{compactc}
There exists reasonably growing function $L=L({1\over r_0})$ such that for a smooth isotopy $K_t$ between knots $k_1$ and $k_2$ in the standard round $S^4$ satisfying $r(k_1)\geq r_0, r(k_2)\geq r_0$, we can obtain a smooth isotopy $K'$ between $k_1$ and $k_2$ with ${1 \over r(K'_t)} \leq L \cdot C(K_t)$ and ${1 \over L \cdot C(K_t)} \leq \Area(K'_t) \leq L \cdot C(K_t)$.
\end{Lem}

\begin{proof}
The only way $1 \over r$ is much larger than $C$ is if the area of $K_t$ is very small, but then, we claim, its diameter is small and therefore the knot is located inside a small ball in $S^4$ that we can scale up like one can do in $R^4$, thus decreasing $1 \over r$. We need to estimate the diameter of $K_t$ in terms of its area. We can do it through estimating the intrinsic diameter of $K_t$ as follows.

Let $\gamma$ be a geodesic in $K_t$ of length $D$, the diameter of $K_t$. We place $D \over r$ points on $\gamma$, each at the distance $r$ from its neighbours. We claim that geodesic balls (in the ambient $S^4$) centered at these points of radius $r \over 10$ do not intersect, we call them $B_i$. This claim follows from the fact that projecting on $K_t$ increases lengths by at most $10 \over 9$ (Lemma \ref{radial}), while $({r \over 10 }+ {r\over 10})({10 \over 9}) < {r}$.

Now we want to bound the volume of $r \over 10$-neighbourhood of $K_t$ in terms of the area of $K_t$. Problem 5-2 from \cite{tubes} gives us the equation $\Vol(N_{r \over 10}(K_t)) = \int_{0}^{0.1r} \const \sin x (\Area(1 -{3 \over 2} \sin^2 x) + 2 \pi \sin^2 x)dx$, where $\Area = \Area (K_t)$ and $r = r(K_t)$. For small $r$ we obtain $\Vol(N_{r\over 10}(K_t)) \leq \const_1 \Area \cdot r^2 + \const_2 \cdot r^4$. We can find a lower bound for the volume of $N{r\over 10}(K_t)$ by adding up the volumes of $B_i$ which yields the inequality $\const_3 {D \over r} \cdot r^4 \leq \const_1 \Area \cdot r^2 + \const_2 \cdot r^4$, or $D \leq r (\const_4 C + \const_5)$, where $C = C(K_t) = {\Area \over r^2}$. We scale up $K_t$ up to some constant intrinsic diameter $D_0$ if its intrinsic
diameter is less than $D_0$ and call the resulting isotopy $K'_t$. (Here the center of the rescaling should not be too close
to $K_t$.) The complexity $C(K'_t)$ does not change by much: we choose $D_0$ small enough so that $ 0.001 C < C' = C(K'_t) < 1000 C$. Therefore, we obtain ${1 \over r'} \leq {\const_4 C' + \const_5 \over D_0} \leq 1000{\const_4 C + \const_5 \over D_0}\leq {\const_6\ C\over D_0}$ as $C$ is uniformly bounded below by a positive constant.

\par

Since $r'$ can not be greater than $\pi$, $\Area' \leq \pi^2 C'$. Also, $\Area' = (r')^2 C' \geq {D_0^2\over \const_6^2 C^2} 0.001 C = {\const_7 D_0^2\over C}.$

However, we do not want to change $k_1$ and $k_2$ during the described
procedure. For this purpose we first observe that the above argument works without any rescaling of $k_1$ or $k_2$ 
(or, more formally, with the identity rescaling of these knots), when
$r_0\geq {D_0\over 2}$ as the intrinsic diameter of a knot $k$ is greater than or equal to its extrinsic diameter $\geq 2r(k).$ Then we notice that if $r_0<{D_0\over 2}$, we can simply replace $D_0$ by a new value, namely, by $D_0(r_0)=2r_0$ and observe that the above argument is still valid without any changes.

\end{proof}

\begin{Lem} \label{step}
Let $k_1, k_2$ be two knots in $S^4$. Let $F \geq \max \{{1 \over r(k_1)}, {1 \over r(k_2)}\}$, and assume $S^4 \setminus N_{1\over 10F}(k_i)$ are ${1\over 10000F}$-close to each other in Gromov--Hausdorff metric, then $\Fl(k_1) < 10^6 \cdot \Fl(k_2)$.
\end{Lem}

\begin{proof}
First we want to show that we can contract small, of length $l \leq {1\over 100 F}$, curves in $C_i = S^4 \setminus N_{1\over 10F}(k_i)$ through curves of length at most $5l$. Let $\gamma$ be such a curve. If $\gamma$ is not in the $l$-neighborhood of $\partial C_i$, then it is contained in a convex metric ball of radius $l\over 2$ in $S^4$ and can be contracted within this ball without any length increase. Otherwise, it can be homotoped along outer normals to $\partial N_{1\over 10F}(k_i)$ to the outer boundary of the $2l$-neighborhood of $\partial N_{1\over 10F}(k_i)$  (or, equivalently, $({1\over 10F}+2l)$-neighborhood of $k_i$). Note that the length of $\gamma$ under the projection increases less than by a factor of $2$ (the second part of Lemma \ref{radial}). Therefore, we can place $\gamma$ in a Euclidean ball of radius $l$ and contract without length increase.

Now we can adapt the argument from section 2 of \cite{nabutovsky3} to prove the lemma. The idea is that in order to contract a contractible closed curve $\gamma_1$ in $C_1$, one can transfer the curve to $C_2$, which is close to $C_1$ in Gromov--Hausdorff metric, without significant length increase, contract the resulting curve $\gamma_2$ there, discretize the contracting homotopy, transfer it back to $C_1$ and ``fill" the discretized homotopy. We transfer curves in the following way: choose points on the curve distance $d$ apart, where $d$ is the Gromov--Hausdorff distance. Then we transfer them to the other space which makes them at most $3d$ apart and connect by a minimizing geodesics. We transfer homotopies similarly and then need to contract geodesic squares of length at most $12d$ to restore the homotopy on the other side.

In order for this program to work one needs to be able to contract without significant length increase ``short" closed curves. Here ``short" means of length that does not exceed the distance between $C_1$ and $C_2$ times an appropriate constant. (It is not difficult to see that one can choose this constant to be equal to  $12$.) In order to see that $\gamma_2$ is contractible in $C_2$ we can construct a contracting homotopy by similarly ``transferring" a homotopy contracting $\gamma_1$ in $C_1$. We refer the reader to the proof of Lemma 2 in \cite{nabutovsky3} for detailed descriptions of such transfers. This argument works, for example, if the distance between  $\partial N_F(k_i),\ i=1, 2$ does not exceed ${1\over 10000F}$.
\end{proof}

\begin{Lem} \label{precompact}
There exists a reasonably growing function $f$ with the following property. Let $S(F,G,L,D) = \{ S^4 \setminus N(k) | r(k) \geq {1 \over F}, L \leq \Area(k) \leq G , \Diam(k) \leq D \}$ the set of knot complements endowed with Gromov--Hausdorff metric, where $N(k)$ denotes the ${r(k)\over 10}$-neighbourhood of $k$. Then there exists a ${1\over 10000F}$-net of size at most $f(F+G+D+{1 \over L})$ on $S(F,G, L, D)$.
\end{Lem}

\begin{proof}
First we are going to construct a ${1 \over 10000F}$-net on $M = (S^4 \setminus  N(k)) \in S(F,G,L,D)$ and obtain a uniform reasonably growing upper bound for sizes of nets on $M$ in the following way.

Let $M' = (S^4 \setminus N_{\epsilon}) \subset M$, where $\epsilon = {1 \over 10}r(k)+{1\over 20000F}$. We can find a finite set of points, call it $T_1$, in $M$ whose ${1\over 100000F}$-neighbourhood contains $M'$. This set will be a subset of a ${1 \over 200000F}$-net on $S^4$, which can be obviously chosen of reasonably growing size as a function of $F$. The subset consists of those points which lie in $M$. The distance between two points in $S^4$ of length less than ${1 \over 100000F}$ doesn't increase by more than the factor of $2$ when regarded as intrinsic distance in $M$ (by the same argument as in the previous lemma), therefore ${1\over 100000F}$-neighbourhood of $T_1$ in the intrinsic metric of $M$ contains $M'$.

Then we find a ${1\over 20000F}$-net of reasonably growing size for $\partial M$ with the induced Riemannian metric, call this net $T_2$. We can do that because we have bounds on curvature, diameter and volume of $\partial M$. If we take the ${1\over 10000F}$-neighbourhood of $T_2$ in $M$ it will contain $M \setminus M'$ due to the triangle inequality. Therefore, $T_1 \bigcup T_2$ is a ${1\over 10000F}$-net for $M$.

Then one can use the standard argument, namely the above mentioned Gromov precompactness criterion, to prove the lemma.
\end{proof}

Let $K=\Sigma\setminus S$ denote the complement
of a $2$-sphere $S$ in a Riemannian $4$-sphere $\Sigma$.
As $K$ is not compact, Lemma 4.2 does not imply that $\Fl(K)<\infty$.  Nevertheless, this is 
so, and will be seen from the following proposition.
Assume, that for some $x$
$r(S)\geq \frac{1}{x}$, the absolute value of sectional curvature of $S$ is less than $1$
(this is always achievable by a rescaling), 
volume of $S$ is $\frac{1}{x^2}$ and diameter of $S$ is $\leq x$. The Cheeger inequality now implies that
the injectivity radius of $S$ is greater than $\const\ \frac{1}{x^2\exp(x)}.$ Let $N$ denotes
the open tubular neighbourhood of $S$ of radius $\frac{r}{C}$ for some constant $C$, and $K_1$ denotes
$\Sigma\setminus N$. Our next proposition
tells us that $\Fl(K)$ and $\Fl(K_1)$ differ by
a factor bounded by a reasonable growing
function of $x$. In particular, as $\Fl(K_1)$
is finite by virtue of Lemma 4.2, $\Fl(K)$
is finite as well.

\begin{Lem} \label{complement} Let $x$, $K$, and $K_1$ be as above.

 There exists a reasonably growing
function $f(x)$ such that 
$$ \Fl(K)\leq f(x)\Fl(K_1).$$
\end{Lem}

\begin{proof}
We need to estimate $\Fl(K)$ in
terms of $\Fl(K_1).$
We would like to prove that if $\gamma$ is
a curve in $K$, then there exists a closed
curve $\gamma_1$ in $K_1$ such that
$\filllength(\gamma_1)$ in $K_1$ will be less than
$f(x) \Fl(K_1) \length(\gamma)$ in $K$ for an appropriate
reasonably growing function $f$. Moreover,
$\gamma$ and $\gamma_1$ can be connected by
a homotopy in $K$ that passes through curves
of length $\leq f(x) \Fl(K_1)\length(\gamma)$.

The main problem is the potential existence of uncontrollably many
very short arcs that loop around $S$ in different directions. As $\gamma$
is contractible in the complement of $S$, the linking
number of $\gamma$ and $S$ is zero. This gives us the idea of ``sliding" these
short loops along $S$ so that they become close to each other, and the great
majority of them can be paired with oppositely directed loops and cancelled.
Then the resulting curve could be pulled away from $S$ without increasing
its length too much. Here is a description of  how this can be accomplished.

Let $i$ denote the minimum of the convexity radii of $S$ and $\Sigma$. Consider a smooth triangulation of $S$
into ``fat" triangles $T_k$ of diameter $i/10$ such that the maximal radius of a metric ball inside
$T_k$ is greater than
$i/g(x)$ for a reasonably growing $g$. (The possibility of having ``large" radii of metric balls in $T_k$ easily follows from the fact that we control the Lipschitz constants
of the exponential map $\exp_p$ and
its inverse at each point $p\in S$ at distances $<i/2$.) Consider ``slabs" $S_k$ made of points of $\Sigma$ 
on normals to $S$ of length not more than $r/C$ starting at a point of $T_k$. We are going to consider all maximal
connected arcs of $\gamma$ that are in some slabs $S_i$. Order them in the same order they appear along $\gamma$ starting from some point of $\gamma$: $S_{j_1}, S_{j_2}, \ldots , S_{j_N}$, and then $S_{j_1}$ again.  Of course,
this list can include multiple copies of the same slab. The curve $\gamma$ might also contain arcs $\beta_k$ after exiting slab $S_{j_k}$ and before entering the next slab.
Denote the maximal (connected) arc of $\gamma$
in $S_{j_k}$ by $\alpha_k$. We need a resonably growing in $x$ upper
bound for the number $N$ of (possibly repeating) slabs in this sequence. Such an
upper bound would be automatic if the length
of each $\alpha_k$ were at least ${i\over 1000}$. To reduce the general case to this situation observe that if $\alpha_k$ has length $<{i\over 1000}$, then it lies close
to the boundary of $S_k$ and can be pushed to
the boundary (and, thus, outside of the open slab
$S_k$) by a radial projection from the center
$p$ of the slab. (As we have a nice control
over Lipschitz constants of $\exp_p$ and $\exp_p^{-1}$, we can think about path homotopies
of arcs in a slab almost as if the slab was
in the Euclidean space.) After pushing out all short
arcs outside of the slabs we will be left
only with long arcs in a controlled number of slabs as desired.

Now our strategy will be to homotop $\alpha_k$ to the minimal geodesic
in $\Sigma$ between the endpoints of $\alpha_k$, and then to push it out of the slab
(leaving possibly only one or two vertical ``whiskers" going up along the boundary. This
happens only if $\gamma$ enters $S_{j_k}$ from
its side. This can happen only if it enters
immediately from $S_{j_{k-1}}$ with trivial
intermediate arc $\beta_{k-1}$. In this case
pairs of
whiskers emanating from the same point $p$ of $\gamma$ and corresponding to adjacent slabs will afterwards be cancelled. Of course, to make this possible  one needs
to go ``up" in the identical way
for both occurrences of $p$ at the boundary of
the adjacent arcs of $\gamma$, which will not present any difficulties.) However, this
homotopy might not be possible if $\alpha_k$ wraps about
$S$ a non-zero number of times. Moreover,
we do not even have control over the winding
number $w_k$ of $\alpha_k$ about $S$ as $\alpha_k$ can come to $S$ arbitrarily close. (To formally define $w_k$ we close $\alpha_k$ by adding the minimizing geodesic in $\Sigma$
between its endpoints.) Note that $N_0=\sum_kw_k$ need not be zero, but $ N_0$
is the winding number of the rest of $\gamma$
about $S$. As the rest of $\gamma$ lies ``far" from $S$ one can majorize $\vert N_0\vert$ 
by a reasonably growing function of $x$.

For each arc $\alpha_k$ denote its point closest to $S$ by $y_k$. Without any loss
of generality we can assume that $y_1$ is the closest to $S$ among all points $y_k$. Now the
idea is to create $\vert N_0\vert+\sum_{k=1}^N\vert w_k\vert$
loops going in positive and negative directions (in accordance to the signs of $N_0$ and $w_k$)  from $y_1$ about $S$ along the same metric circle
of radius $dist(y_1, S)$ centered at the point of $S$ nearest to $y_1$. The total length, however, will be less than the length of $\gamma$.

The next step will be to move $\sum_{i=2}^N\vert w_i\vert$ circles currently
attached to $\alpha$ at $y_1$ to $y_2,\ldots y_N$. We will be moving 
these circles one by one, first, to $y_2$, then to $y_3$, etc. until they reach the desired destinations. Moving a circle from $y_i$ to
$y_{i+1}$ is equivalent to contracting (in $K$)
of the closed curve $\rho$ formed by two copies of
the arc of $\alpha$ between $y_i$ and $y_{i+1}$
traversed in opposite directions separated by
oppositely oriented circles around $S$ at $y_i$ and $y_{i+1}$. The winding number of this closed curve is $0$, and therefore it is null-homotopic in $K$. 
We start from homotoping two metric circles on both ends of the arc along metric circles based at points of the arc.
We are moving the base points along the arc until they exit the corresponding slabs, and the considered curve $\rho_1$
that needs to be contracted will be outside of the ${r\over C}$-neighbourhood of $S$, that is in $K_1$.
The curve $\rho_1$ gets shorter as we consider shorter arcs of $\alpha$ but can increase its length because we replace two very short
metric circles at both ends of the curve by metric circles of radius slightly larger ${r\over C}$. However,
as the segments of $\rho$ coming from $\alpha$ are sufficiently long, this will not increase the length by more than
a $\const$ factor. Once the considered closed curve is in $K_1$, its contractibility will be controlled by $\Fl(K_1)$.

As we
move these small loops one by one, these temporary length
increases do no accumulate in the process.

Once we moved these loops in the desired positions, we are going to have the new arcs $\tilde\alpha_k$ that
look like $\alpha_k$ plus $\vert w_k\vert$ small circles oriented so that $\tilde\alpha_k\bigcup$ the geodesic segment connecting its endpoints has winding number
$0$ and is null-homotopic. The key idea now is that the complement to a $2$-plane in the Euclidean $4$-space has $\Fl$ equal to $1$.
As slabs are nicely bi-Lipschitz homeomorphic
to convex domains in $\mathbb{R}^4\setminus\mathbb{R}^2$, contractible 
closed curves in slabs can be contracted with
only a nicely controlled length increase. 
This means that we can now implement our original program - at least for slabs $2,3,\ldots , N$ - and to push all arcs
of $\alpha$ but $\tilde\alpha_1$ outside of $N$
into $K_1$.

We are facing an additional problem with
$\tilde\alpha_1$. At the end of the reduction process
we are going to get an arc with $\vert N_0\vert$ small
circles about $S$ at $y_1$ which obviously has winding number $N_0$.
If not these circles, we would homotop
this arc first to the minimizing geodesic, then
to an arc $\bar\alpha_1$ consisting of at most two ``whiskers"
in $\partial S_{j_1}$ leading from the endpoints of $\alpha_1$ to the boundary of $K_1$ and a short arc in the boundary of $K_1$.
(The whiskers would be later cancelled with the same but oppositely oriented whiskers coming from the adjacent slabs.) The idea is to homotop the arc $\tilde\alpha_1$ with $\vert N_0\vert$ small circles around
$S$ to $\bar\alpha_1$ together with $\vert N_0\vert$ (much larger)
circles around $S$ of radius $r\over C$ that are attached to a middle point of $\bar\alpha_1$.
(As $\vert N_0\vert$ is nicely controlled in terms of $x$, the resulting arc in $K_1$ will have length $\leq h(x)\length(\gamma)$ for a reasonably growing function $h$.)  For this purpose
we first ``forget" about $\vert N_0\vert$ small circles,
and homotop $\tilde\alpha_1$ to $\bar\alpha_1$
with $\vert N_0\vert$ small circles still attached to
a point $y$ in its middle part via two oppositely oriented copies
of a minimizing geodesic connecting $y$ and $y_1$. Then we start rescaling these circles
to a larger and larger radii simultaneously cancelling larger and larger arc of the geodesic $yy_1$ starting at $y_1$. At the end
we are going to get the desired arc that
after the subsequent cancellation of the two whiskers at its ends will be in
$K_1$. These completes the demonstration
of the upper bound for $\Fl(K)$ in terms
of $\Fl(K_1).$

\end{proof}

Now, we are going to observe that:

\begin{Pro} \label{standard}
The functions $\Fl(M(Q_n)\setminus S_n)$ and $\Fl(S^4\setminus \phi_n(S_n))$
are similarly growing functions of $n$.
\end{Pro}

\begin{proof}
This assertion immediately follows
from Lemma \ref{knots} and the obvious observation that a $L$-bi-Lipshcitz homeomorphism changes $\Fl$ by at most the factor of $L^2$. 
\end{proof}

Now it remains just to prove the following proposition that will complete our proof of Theorem \ref{smoothmain}:

\begin{Pro} \label{dehn}
The function $\Fl(M(Q_n)\setminus S_n)$ is not
reasonably growing.
\end{Pro}

\begin{proof} 
We are going to prove this theorem by contradiction. We are going to
assume that $\Fl(M(Q_n)\setminus S_n)$ is reasonably growing.

In view of Theorem 2.1 we will be proving that the maximal length of closed curves in
a nearly optimal homotopy contracting the closed curve corresponding to the generator
$x_{n+1}$ grows similarly to the number of 
%
elements in a sequence of words in $P_n$ starting from $x_{n+1}$ and ending at the trivial word,
such that the maximal length of these words where we do not count $t$s and $t^{-1}$s
grows as the maximal
length of closed curves during the contracting homotopy.

The Riemannian manifold $M(Q_n)$ was
built starting from the presentation  $2$-complex for the finite presentation
$Q_n$ of the trivial group. Geometrically,
one can think about it as about the (smoothed-out) boundary
of a tubular neighbourhood of the presentation complex sitting in $\mathbb{R}^5$, where
all $2$-cells have reasonably simple geometries. Formally speaking, this means
that they are bi-Lipschitz homeomorphic to
the corresponding Euclidean discs with Lipschitz constants being slowly (in particular, reasonably) growing functions of $n$. Geometrically, the difference between
the tubular neighbourhood and its boundary consists in deletion of cells of dimensions $3,\ 4$ and $5$ - all with simple geometry.
Therefore, the standard argument from a graduate course in topology that adding/removing cells of dimension $\geq 3$ does not change
the fundamental group, can be easily adapted
to demonstrate that each null-homotopy
of a contractible closed curve in the boundary
inside the tubular neighbourhood could be easily projected to the boundary of the tubular neighbourhood without increasing lengths of curves by much. Further, the existence of an obvious deformation retraction of the 
tubular neighbourhood to the presentation
complex itself implies that the filling constant 
$\Fl$ of the boundary of the tubular neighbourhood will be ``almost the same" as
$\Fl$ of the tubular neighbourhood , which
will be ``almost the same" and $\Fl$ of the
presentation complex. This is a very general
argument that pertains not only to $Q_n$ but
to all finite presentations of (trivial
or non-trivial) groups with ``small" 
presentation length. 

This argument would pertain to
$M(Q_n)$. Yet observe that we are interested
in $M(Q_n)\setminus S_n$. Geometrically, we
delete the center from one of $2$-cells of
the presentation complex of $Q_n$, the
$3$-dimensional ``equatorial" disc from
the corresponding $5$-dimensional cell of
the tubular neighbourhood, and the $2$-sphere
from $S^4$ in the boundary of the tubular
neighbourhood which is the boundary of the
aforementioned $5$-cell. As these pairs ($5$-cell, $3$-cell) and $(S^4,S^2)$ have 
very simple geometry, the above argument
implies that $\Fl(M(Q_n)\setminus S_n)$ grows
similarly to $\Fl$ of the presentation
complex of $Q_n$ with the puncture at the center of the disc corresponding to the relator $t$. Here we can assume without any loss of generality (an up to a reasonably growing factor) that the metric of the presentation complex is the path metric, where all generators are circles of length $1$, and
all $2$-cells before their attachment are Euclidean discs with the length of the circumference equal to the length of the corresponding relator.

Now we are going to try to relate $\Fl$ for
such complexes with the Dehn function and the filling length function
of the corresponding group presentation minus the relator corresponding to the punctured $2$-disc. In our case we are
talking about $Q_n$ minus the relator $t$, 
i.e. $P_n$. 

Consider first the presentation 
$2$-complex for $Q_n$ without any $2$-cells containing the punctures, that is, the presentation complex for $P_n$. 
It is well-known that the filling length in a presentation complex 
of a finitely presented group is closely related with
the Dehn function of the finite presentation
- see, for example, \cite{GR} for details. In particular, one of these two functions
is reasonable if and only if the other is.

Our situation is somewhat different as rather than studying different
words in a fixed finite presentation we study the situation
when the set of finite presentations is variable but one considers only one word in each finite presentation.
Nevertheless, we can apply similar ideas.
Consider a nearly optimal homotopy $\{\gamma_t\}_{t\in [0,1]}$ in the presentation complex 
contracting the curve $\gamma=\gamma_0$ in the $1$-skeleton corresponding to $x_{n+1}$
to a point. By a ``nearly optimal homotopy"
we mean a homotopy, where the maximal length of closed curves during the
homotopy is almost minimal possible. Triangulating all $2$-cells as
the cones with a new vertex at the centre of the cell over a triangulation
of the boundary, and
replacing this homotopy by its simplicial approximation, we can assume
that it is simplicial. We can think about this homotopy as about a sequence
of simplicial closed curves in the $1$-skeleton of the triangulation
of the presentation complex, where the symmetric difference of 
each pair of consecutive closed curves is the boundary of a $2$-simplex
of the triangulation. Assuming that the triangulation of the presentation
complex has the minimal number of $2$-simplices, this number of
$2$-simplices is proportional to the length of $P_n$. Therefore, the
number of all distinct simplicial curves of length bounded by a reasonably
growing function of $n$ is reasonably growing. An optimal
homotopy does not revisit the same closed curve twice. So, the curves
in this homotopy pass through the center at a $2$-cell 
reasonably many times. Without any loss of generality no closed
curve in the homotopy passes through more than one center of a $2$-cell.

Now forget about the chosen triangulation.
At each moment of the homotopy when the closed curve does not pass
through the center of a $2$-cell, it can be radially projected from
the centers of $2$-cells to a closed curve in the $1$-skeleton of
the presentation complex. Each such curve corresponds to 
a word in the group up to a conjugacy corresponding to a choice
of the shortest path connecting the base point and the initial
point of the curve. The projected curves change continuously until
homotopy passes through the center of a cell. 
The corresponding words change also only in those moments, and each change
can be described as the deletion/insertion of a conjugate to the
relator corresponding to the $2$-cell. This would mean that there exists 
a van Kampen diagram with $x_{n+1}$ on the boundary with a reasonably
growing as a function of $n$ number of cells that yields a contradiction
with Theorem 2.1.

Yet we also need to account for a presence of the $2$-cell
corresponding to the relation $t$ with the
puncture at its center. 
So, consider again the nearly optimal homotopy the closed 
simplicial curve $\gamma$ representing $x_{n+1}$ to a point
in the punctured metric presentation $2$-complex
for $Q_n$. At each moment $t$ of the homotopy we consider all arcs of $\gamma_t$
inside this cell. Without any loss of generality we can assume that the set of such arcs
is finite. Also, without any loss of generality we can assume that none of closed curves $\gamma_t$ is contained inside the punctured $2$-cell.

%

We can continuously push the arcs of this homotopy in the
punctured cell to its boundary by the radial projection
with the center at the puncture. As a result, we can get
an uncontrollably large number of letters $t$ and its inverse.
The number of $x$s, however,
remains reasonably growing as a function of $n$. 

Now the application of Lemma \ref{notee} 
completes the proof of the proposition.

\end{proof}

{\bf Remark.} In this remark we would like to sketch a different version of the proof of Theorem 1.1.
Its detailed exposition would be cumbersome and tedious. Yet it has a simpler logical
structure, does not use Lemma 4.8, Lemma 2.2, as well as Karcher's smoothing with mollifiers on Riemannian manifolds in the proof of Lemma 3.7. 

The basic idea now is that we control $C^1$ norm of the curvature of manifolds $M_t$ when they are constructed, and make the isotopy connecting them continuous in $C^2$-topology.
(In fact, one can even construct $M_t$ as hypersurfaces in $\mathbb{R}^5$ and ``code" them by functions $f_t$ that vanish on $M_t$, change from $-1$ to $1$ in the tubular neighbourhood of $M_t$ of radius $r(M_t)/2$, and are equal to $-1$ or $1$ outside this neighbourhood - compare \cite{nabutovsky2}. This approach would eliminate a significant amount of Riemannian geometry from the proof.) This isotopy can be discretized, so that the diffeomorphisms $\phi_{t_i, t_{i+1}}$ between consecutive manifolds in the discretization , the differentials of
$\phi_{t_i,t_{i+1}}$, as well as their inverses, and the differentials of their inverses have
Lipschitz constants close to $1$. As the corollary, such a diffeomorphism will map
the boundary of the tubular neighbourhood of the $2$-knot in $M_{t_i}$ to
a normal variation of the boundary of the tubular neighbourhood of the $2$-knot in
$M_{t_{i+1}}$, so that there exists an explicit and well controlled isotopy in $M_{t_{i+1}}$
connecting these two hypersurfaces in $M_{t_{i+1}}$. This isotopy can be used to alter $\phi_{t_i,t_{i+1}}$ in a controlled way so that it bijectively maps the tubular neighbourhood of  the $2$-knot in $M_{t_i}$ to the tubular neighbourhood of the $2$-knot in 
$M_{t_{i+1}}$. The advantage is that now we have an analog of Proposition 4.9 for complements
of tubular neighbourhoods of $2$-knots, and do not need to use $Fl$ for complements of
the $2$-knots at all, but only $Fl$ for the complements of the tubular neighbourhoods of the $2$-knots. This will eliminate
the need for Lemma 4.8, and, as there will be no difficulty with the punctured $2$-cell
at the end of the proof of Proposition 4.10 eliminating the need for Lemma 2.2.
(The punctured $2$-cell will lack a disc of a definite size, so the radial projection
will increase the length by a reasonably growing factor.)

We will need to demonstrate that the number of $C^2$-close  manifolds $M_{t_i}$ in the discretization of a path continuous in the $C^2$-topology is bounded by a reasonably growing
function of $n$. This will follow from the proof of the relevant version of the Ascoli-Arzela
theorem once we have reasonably growing bounds for the relevant parameters of manifolds $M_t$. To achieve this one would need to have a look at the isotopies  of smooth embedded curves
corresponding to the Andrews-Curtis operations of type 1 and 3 (see section 3.2), and the surgeries creating the manifolds $M_t$. Whenever it is obvious that all this can be straightforwardly accomplished,
we do not want to discuss the tedious details.

\section{PL case}
In this section we will abuse notation by using the same symbol for both a simplicial complex and its geometric realization. For example, $S^4$ will denote both the space and the simplicial complex $\partial \Delta^5$. By ``subdivision'' (or refinement) we will always mean rectilinear subdivision. Let $L \subset K$ be simplicial complexes. By $N(L,K)$ we will denote the simplicial neighbourhood of $L$, and by $C(L,K)$ the simplicial complement. A subdivision $K'$ of $K$ we will call derived near $L$, if it is obtained from $K$ by selecting a new vertex in every simplex that neither lie in $L$, nor in $C(L,K)$. We start by defining combinatorial filling length similarly to the Riemannian case: 
\begin{Def}
For a triangulated compact space $M$ and a simplicial closed curve $\gamma$, let $H(\gamma)$ be the set of all simplicial homotopies contracting $\gamma$ (a simplicial homotopy is a finite ordered collection of simplicial curves, where two neighbouring curves differ by a $2$-simplex). For $h \in H(\gamma)$ let $\fll(h)$ denote maximal length of a closed curve in $h$ (length of a curve here is the number of $1$-simpleces). Define $\filllength(\gamma) = \inf_{h \in H(\gamma)} \fll(h)$. Finally, let $\Fl(M)$ be the supremum of $\filllength(\gamma) \over \length(\gamma)$ over all closed contractible curves.
\end{Def}

Alternatively to $\Fl$ we could use the concept of visible presentations of the fundamental group defined in \cite{ln3}, but we decided to go with $\Fl$ to make this section similar to the previous one. We will also need the following facts relating PL topology to combinatorics.

\begin{Lem} \label{facts}
Let $[K]$ denote the number of simplices in a complex $K$.
\begin{enumerate}
    \item Given a linear embedding $K_1 \to K$ of triangulated spaces, we can refine the triangulation on $K$ to contain $K_1$ as a full subcomplex so that the number of simplices is a reasonably growing function of $\max ([K], [K_1])$.
    \item Given two refinements $K_1, K_2$ of the same triangulation $K$, we can find a common refinement with the number of simplices bounded by $d[K_1][K_2]$, where $d$ is a constant which only depends on the dimension.
    \item Let $K_1, K_2$ be two triangulations of $(M,N)$, a pair of closed manifolds, $N \subset M$. Let $L_1, L_2$ be the restriction of those triangulations to $N$. Let $K'_i$ be a derived subdivision near $L_i$. Assume that $L_i$ is a full subcomplex of $K'_i$, $i=1,2$. Then there is a PL-homeomorphism $K_1 \to K_2$ fixing $N$ and sending $N(L_1, K'_1)$ onto $N(L_2, K'_2)$ of complexity at most $d[K]^4$, where $K$ is a common refinement of $K_1, K_2$ and $d$ is some constant which depends only on dimension of $M$.
    \item There exists a reasonably growing function $f$ with the following property. Let $K'$ be a subdivision of $K$, a triangulation of a compact manifold. Then if the manifold has no boundary, $K$ and $K'$ are related by a sequence of bi-stellar moves of length at most $f([K'])$. If boundary is nonempty, then the same result holds for shellings/inverse shellings instead of bi-stellar moves.
\end{enumerate}
\end{Lem}

\begin{proof}
\begin{enumerate}
    \item Assuming $K_1$ has one simplex we can easily refine the triangulation $K$ to include it. So the result follows by induction on the number of simplices in $K_1$. We also derive the triangulation once to make sure $K_1$ becomes a full subcomplex. 
    \item Take a simplex $v_1$ in $K_1$ and a simplex $v_2$ in $K_2$. Their intersection is a convex polyhedron because they both lie in some larger simplex of $K$. Define $d$ to be the maximal number of simplices we need to triangulate an intersection of two rectilinear simplices in the given dimension.
    \item That's an effective version of Theorem 3.8 from \cite{pl} and we will follow the proof given in that reference. We are going to show that there is an effective homeomorphism $K_1 \to K$ such that $N(L_1, K'_1)$ is mapped to $N(L, K')$, where $L$ is the restriction of $K$ to $N$ and $K'$ is a derived subdivision near $L$ , and that would complete the proof by composition. We define $h: K_1 \to [0,1]$ by sending vertices of $L_1$ to $0$, all other vertices to $1$ and extending linearly. For example, $h^{-1}(1) = C(L_1, K_1)$. We can find a small enough $\epsilon$ such that $h^{-1}([0,\epsilon])$ has no vertices of $K \setminus L$. We can build $K'$ and $K'_1$ by placing the new vertices on $h^{-1}(\epsilon)$ (note that any two derived subdivisions for the same triangulation are combinatorially equivalent). Then the geometric realizations of $N(L_1, K'_1)$ and $N(L, K')$ are the same, and coincide with $h^{-1}([0,1])$. Since $K'$ and $K'_1$ are subdivisions of the same triangulation $K_1$ we can use part 2 of this lemma to obtain that complexity of $N(L_1, K'_1) \to N(L, K')$ is bounded by $d'[K]^2$, for some constant $d'$ depending only on dimension. And taking composition with the similarly found homeomorphishm $N(L_2, K'_2) \to N(L, K')$ we get the complexity bound $d[K]^4$ (applying part 2 again), for some constant $d$ depending only on dimension.

    \item The result stated without the bound on the number of moves is due to Pachner \cite{pachner}, but the proof, see a nice exposition by Lickorish \cite{lickorish}, is constructive so one can get an effective (reasonably growing) estimate.
\end{enumerate}
\end{proof}

\par\noindent
{\bf Remark.} The standard way to ensure that $L_1,\ L_2$ are full subcomplexes of $K_1, K_2$
in Lemma 5.2.3 is to replace $K_1, K_2$ by their derived subdivisions before applying this lemma.
This would increase the number of simplices in $K_1, K_2$ by not more than a factor that  depends only on the dimension. Therefore, below we will automatically apply Lemma 5.2.3 to the derived sudivisions of the relevant triangulations rather than to the triangulations themselves.

\begin{Lem} \label{plfill}
There exists a reasonably growing function $f$, such that for a PL-homeomorphism $K \to L$ between compact manifolds of complexity less than $c$, $\Fl(L) < f(\Fl(K)+c)$.
\end{Lem}

\begin{proof}
By Lemma \ref{facts}.4 there exists a sequence of elementary transformations (bi-stellar or shellings/inverse shellings) of reasonable length between $K$ and $L$. Each elementary transformation changes $\Fl$ by at most a constant number independent of $K$ or $L$. The lemma follows.
\end{proof}

\begin{Lem} \label{composition}
Let $f:K \to L$, $g:L \to N$, be PL homeomorphisms between manifolds with fixed triangulations. Then $C(g \circ f) \leq d C(f) C(g)$, where $d$ is a constant depending only on dimension.
\end{Lem}

\begin{proof}
Let $L_1$ be a subdivision of $L$ of size $C(g)$ such that $g$ is linear on every simplex. Let $L_2$ be the subdivision of $L$ of size $C(f)$ coming from $K$, i.e. we subdivide $K$ to make $f$ linear and push this triangulation forward. Then we can take a common subdivision of $L_1$ and $L_2$ of size bounded by $d[L_1][L_2] = dC(g)C(f)$ from Lemma \ref{facts}.2 and pull it back on $K$ to finish the proof.
\end{proof}

We are going to triangulate the sequence of knots constructed previously (see sections 3.2, 3.3) and use it to prove the main result. We summarize the results we need about this sequence in the following proposition which is a PL analogue of Proposition \ref{knots} combined together with
a PL-version of analogs of Proposition 4.9 and Proposition 4.10 for the complements to the tubular
neighbourhoods of $2$-knots (and not for
the complements to $2$-knots themselves as in section 4).

\begin{Pro} \label{plknots}
There exists a sequence of PL-unknots $k_i: S^2 \to S^4$ with the following properties. There exist triangulations $(S^4, k_i)$ with a reasonable number of simplices (as a function of $i$) such that the simplicial neighbourhoods $N_i$ of $k_i$ are regular, and $\Fl(S^4 \setminus N_i)$ is not a reasonably growing function of $i$.
\end{Pro}

\begin{proof}

First we are going to show that we can triangulate $M(Q_n)$ with reasonable number of simplices. (Here and below we are using the notations introduced in sections 3.2, 3.3.) For any presentation we can embed its presentation complex in $\mathbb{R}^5$ and triangulate a small neighbourhood of it with a reasonable number of simplices such that the triangulation, call it $\tilde M(Q_n)$, shells to the presentation complex, denoted $P(Q_n) \subset \tilde M(Q_n)$. The boundary of $\tilde M(Q_n)$ is a triangulation of $M(Q_n)$. We can do the same for any finite presentation.
\par
Note that for any presentation $Q$ all three of $M(Q)$, $\tilde M(Q)$, and $P(Q)$ have ``similarly growing'' $\Fl$; in fact, for either $Q=Q_n$, or $Q=P_n$, the three $\Fl$ are similarly growing functions of $n$ for the following reasons. By a dimensional argument both curves and homotopies of curves in $\tilde M(Q)$ can be pushed off $P(Q)$, and then they can be pushed further to the boundary, $M(Q)$, by inverse shellings.

\par
We can also make sure that the $2$-knot $S_n$ is a subtriangulation and a full subcomplex, derive the triangulation near $S_n$, and we can do the same for $D_n$ -- the obvious $3$-disk filling of $S_n$ in $\tilde M (Q_n)$, that is $D^3 \times \{c \}$ in the last $2$-handle. For the same reason as in the previous paragraph, $\Fl(M(Q_n) \setminus N(S_n))$ and $\Fl(\tilde M(Q_n) \setminus N(D_n))$ are similarly growing.

\par

There is an obvious PL-homeomorphism of small complexity from $\tilde M(Q_n) \setminus N(D_n)$ to $\tilde M(P_n)$ where we retract the last $2$-handle. Therefore, we have a chain of spaces: $P(P_n)$, $\tilde M(P_n)$,  $\tilde M(Q_n) \setminus N(D_n))$, $M(Q_n) \setminus N(S_n)$ with similarly growing $\Fl$, while $\Fl(P(P_n))$ is not reasonably growing by Theorem \ref{groups}. Therefore $Fl(\tilde M(Q_n) \setminus N(D_n))$ is not reasonably growing.

\par

Finally, notice that Andrews-Curtis operations induce obvious PL-homeomorphisms of controllable complexity for the following reason. Take two presentations $R_1, R_2$ related by an Andrews-Curtis operation. The isotopy of the attaching circle relating $\tilde{M}(R_1)$ to $\tilde{M}(R_2)$ fixes all of the circle except for a small interval that traces out a disk (in the case of the handle slide it goes over a $2$-handle, and if it is a free cancellation it stays on the $1$-handle). This disk is broken into triangles by the triangulation of the manifold, where the number of triangles is bounded by the number of simplices. So this isotopy is a sequence of ``embedded Pachner moves'' discussed in the introduction. Each such move is obviously effectively extendable as an ambient isotopy. Therefore the complexity of the homeomorphism $\tilde{M}(R_1) \to \tilde{M}(R_2)$ corresponding to the handle slide depends reasonably on the number of simplices we started with. We can restrict it to the boundary $M(R_1) \to M(R_2)$.

In this way we can construct a sequence of homemorphisms of reasonable complexity between $M(Q_n)=M(Q_n^{(1)})$, $M(Q_n^{(2)})$,...,$M(Q_n^{(m)})=S^4$, where $m\leq\const 2^n$. Then the composition of all of them, called $\phi_n$, is of reasonable complexity by Lemma \ref{composition}. Lemma \ref{facts}.3 gives us a PL-homeomorphism of $S^4$ into itself of reasonable complexity that maps $\phi_n(N(Q_n))$ to $N(\phi_n(Q_n))$ and fixes $Q_n$. Now we can use Lemma \ref{plfill} to finish the argument.
\end{proof}

Now we state and prove a PL analogue of Proposition \ref{comfill}.

\begin{Pro} \label{plcomfill}

Let $k^{(i)}$ be a sequence of 2-unknots in $S^4$ that can be untied through $2$-knots of complexity bounded by a reasonably growing function of $i$ and the isotopies realizing those minimal untyings are $m$-locally extendable, then $\Fl(S^4 \setminus N(k^{(i)}))$ is a reasonably growing function of $i$.
\end{Pro}

\begin{proof}
For $k^{(i)}$ let $K$ (index suppressed) be an $m$-effectively locally extendable isotopy to the trivial embedding of reasonable complexity. We can construct local extensions of $K$: $\tilde K^{(1)}, ..., \tilde K^{(n)}$, for time intervals $[t_0, t_1], [t_1, t_2], ... , [t_{n-1}, t_n]$ with complexities bounded by a reasonably growing function. We don't have control over how big $n$ is.

For each $t_i$ we will subdivide $S^4$ so that a triangulation of $K_{t_i}(S^2)$ with the minimal number of simplices will be a full subcomplex with a reasonable number of simplices (Lemma \ref{facts}.1). Then we take the first derived subdivision of this triangulation so that the simplicial neighbourhood, call it $N_i$, of $K_{t_i}(S^2)$ will be regular. Define $\phi_i = \tilde K^{(i)}_{t_{i+1}}$, it is a PL-homeomorphism of reasonable complexity. Now we want to compare $N'_i := \phi_{i}(N_i)$ to $N_{i+1}$, both regular neighbourhoods of $K_{i+1}$. Let $C'_i, C_{i+1}$ be the complements of these spaces. Then by by Lemma \ref{facts}.3 there exists a PL-homeomorphism of reasonable complexity between $C'_i$ and $C_{i+1}$.
\par
The total number of simplicial complexes of size $s$ is a reasonably growing function of $s$. Therefore, we have a reasonable number of combinatorially distinct spaces $C'_i, C_j$ ``related'' by PL-homeomorphisms of reasonable complexity. ``Related" here means that there is a sequence of PL-homeomorphisms of reasonable complexity between them and the length of the sequence does not exceed the number of combinatorially distinct $C'_i, C_j$. Then by Lemma \ref{composition} there is a PL-homeomorphism of reasonable complexity between any $C'_i$ and $C_j$. Therefore, since $\Fl(C_n)$ is small, we get that $\Fl(C_0)$ is a reasonably growing function by Lemma \ref{plfill}.
\end{proof}

The two propositions above immediately imply Theorem \ref{plmain}. Notice that the only place in the proof of Proposition \ref{plcomfill} where we used the effective local extendability condition was to obtain effective PL-homeomorphisms $\phi_i$. Therefore, we can obtain the following version of the above proposition.

\begin{Pro} \label{plcomfill2}
Let $f$ be a reasonably growing function. Let $k^{(i)}$ be a sequence of 2-knots in $S^4$ for each element of which there exists a sequence of knots $k^{(i)}_1,...,k^{(i)}_{n(i)}$ with the following properties. The first element is $k^{(i)}$, the last is the trivial embedding. For consecutive knots, $(S^4,k^{(i)}_j), (S^4,k^{(i)}_{j+1})$ are related by a PL homeomorphism of complexity bounded by $f(i)$. Finally, $n(i)$ is a reasonably growing function of $i$. Then $\Fl(S^4 \setminus N(k^{(i)}))$ is a reasonably growing function of $i$.
\end{Pro}

 Proposition \ref{plknots} and \ref{plcomfill2} imply Theorem \ref{plcomb}.

{\bf Acknowledgements.} This research has been partially supported by NSERC Accelerator and Discovery Grants of A. Nabutovsky. Research of B. Lishak was also supported by the Australian Research Council's Discovery funding scheme (project number DP160104502). We would like to thank the anonymous referee for numerous comments that helped to improve the exposition.

\bibliographystyle{alpha}

\bibliography{bibliography}

\normalsize

\begin{tabbing}
\hspace*{7.5cm}\=\kill
Boris Lishak                        \>Alexander Nabutovsky\\
School of Mathematics and Statistics           \>Department of Mathematics\\ 
University of Sydney               \>University of Toronto\\
New South Wales,                  \>40 St. George st.,\\
2006             \>Toronto, Ontario M5S 2E4\\
Australia                              \>Canada\\
email: boris.lishak@sydney.edu.au \>alex@math.toronto.edu\\
\end{tabbing}

\end{document}